\documentclass{amsart}
\usepackage{psfig}
\usepackage{amscd}
\usepackage{amsmath}
\usepackage{graphicx}
\usepackage{amsfonts}
\usepackage{amssymb}
\newtheorem{theorem}{Theorem}[section]
\newtheorem{corollary}[theorem]{Corollary}
\newtheorem{lemma}[theorem]{Lemma}
\newtheorem{proposition}[theorem]{Proposition}
\theoremstyle{definition}

\newtheorem{remark}[theorem]{Remark}

\theoremstyle{remark}
\newtheorem*{acknowledgements}{Acknowledgements}

\newlength{\displayboxwidth}
\renewcommand{\theenumi}{\roman{enumi}}

\numberwithin{equation}{section}
\newlength{\qedskip}

\newcommand{\frogcapCorrelationCoefficients}{}

\newcommand{\frogcapReptileGenealogy}{}

\newcommand{\bs}{\boldsymbol{\sigma}}
\newcommand{\br}{\boldsymbol{\rho}}
\newcommand{\bp}{\boldsymbol{\pi}}
\def\openone
{\mathchoice
{\hbox{\upshape \small1\kern-3.3pt\normalsize1}}
{\hbox{\upshape \small1\kern-3.3pt\normalsize1}}
{\hbox{\upshape \tiny1\kern-2.3pt\SMALL1}}
{\hbox{\upshape \Tiny1\kern-2pt\tiny1}}}
\makeatletter
\let\subsubsubsectionname\@empty
\newcounter{subsubsubsection}[subsubsection]

\def\l@subsubsubsection{\@tocline{4}{0pt}{1pc}{9pc}{}}

\def\subsubsubsection{\@startsection{subsubsubsection}{4}%
  \z@{.5\linespacing\@plus.7\linespacing}{-.5em}%
  {\normalfont\itshape}}
\AtBeginDocument{%
  \@for\@tempa:=-1,0,1,2,3,4\do{%
    \@ifundefined{r@tocindent\@tempa}{%
      \@xp\gdef\csname r@tocindent\@tempa\endcsname{0pt}}{}%
  }%
}
\def\@writetocindents{%
  \begingroup
  \@for\@tempa:=-1,0,1,2,3,4\do{%
    \immediate\write\@auxout{%
      \string\newlabel{tocindent\@tempa}{%
        \csname r@tocindent\@tempa\endcsname}}%
  }%
  \endgroup}
\makeatother
\long\def\TeXButton#1#2{#2}%

\begin{document}
\title[Compactly supported wavelets and Cuntz relations]{Compactly
supported wavelets and representations of the Cuntz relations}
\author{Ola~Bratteli}
\address{Department of Mathematics\\
University of Oslo\\
PB 1053 -- Blindern\\
N-0316 Oslo\\
Norway}
\email{bratteli@math.uio.no}
\author{David E. Evans}
\address{School of Mathematics\\
University of Wales, Cardiff\\
PO Box 926\\
Senghenydd Road\\
Cardiff CF2 4YH\\
Wales, U.K.}
\email{EvansDE@cardiff.ac.uk}
\author{Palle~E.~T.~ Jorgensen}
\address{Department of Mathematics\\
The University of Iowa\\
14 MacLean Hall\\
Iowa City, IA 52242-1419\\
U.S.A.}
\email{jorgen@math.uiowa.edu}
\thanks{Research supported by the University of Oslo.}
\begin{abstract}
We study the harmonic analysis of the quadrature mirror filters coming from
multiresolution wavelet analysis of compactly supported wavelets. It is known
that those of these wavelets that come from third order polynomials are
parametrized by the circle, and we compute that the corresponding filters
generate irreducible mutually disjoint representations of of the Cuntz algebra
$\mathcal{O}_{2}$ except at two points on the circle. One of the two
exceptional points corresponds to the Haar wavelet and the other is the unique
point on the circle where the father function defines a tight frame which is
not an orthonormal basis. At these two points the representation decomposes
into two and three mutually disjoint irreducible representations,
respectively, and the two representations at the Haar point are each unitarily
equivalent to one of the three representations at the other singular point.
\end{abstract}
\maketitle

\section{\label{Intro}Introduction}

\setlength{\displayboxwidth}{\textwidth} \addtolength{\displayboxwidth
}{-2\leftmargini} In this paper we show that wavelets may be constructed from representations of
two systems of operator relations, one on $L^{2}\left(  \mathbb{R}\right)  $
and one on $L^{2}\left(  \mathbb{T}\right)  $, for the case of one real
dimension. Focusing on the case of compact support, the analysis reduces to a
certain finite-dimensional matrix problem which is especially amenable to an
algorithmic and computational approach. The associated algorithms are worked
out in detail for a variety of examples which includes the Daubechies wavelet,
and which also reveals some perhaps unexpected symmetries.

One
benefit from the representation theoretic approach to wavelets is
that it provides a coordinate-free way of making precise notions of
irreducibility which occur in the wavelet literature without always having
precise definitions. Specifically, examples in $L^{2}\left(  \mathbb{R}%
^{d}\right)  $, for $d>1$, may occasionally be reduced to simpler examples in
one dimension, i.e., in $L^{2}\left(  \mathbb{R}\right)  $, by a tensor
product construction, but this analysis depends on the chosen spatial
coordinates in $\mathbb{R}^{d}$, while the representation-theoretic approach
in the present paper does not.

One of our results, Corollary \ref{CorPolynomial.3}, specifies in a general
context (for compactly supported wavelets in
$\mathbb{R}^{\nu}$) a decomposition
formula (finite orthogonal sums of irreducible representations) for the
representation associated with a system of high-pass/low-pass filters which
generate the wavelets in question.

It has been known for some time that a class of convolution
operators
from signal analysis,
called 
subband
filters,
satisfy certain operator relations \cite[Lemma 2.1]{Wic93}. Perhaps it is less
well known among experts in
multiresolution
wavelet theory that these operator relations were
introduced in $C^{\ast}$-algebra theory by
J. Dixmier
\cite[Exemple 2.1]{Dix64}
and
J.~Cuntz \cite{Cun77} several decades ago,
and the $C^{\ast}$-algebra they generate is now called the Cuntz algebra of
order $N$ and is denoted by $\mathcal{O}_{N}$,
where $N$ is the scale of the resolution.
This algebra is independent of
the particular scale-$N$ multiresolution wavelet, but the unitary equivalence
class of the corresponding representation may depend on the wavelet. The
detailed structure of these representations has, however, so far only been
worked out in the single case of the Haar wavelet (see below). The purpose of
the present paper is to work out the structure of these representations for
all compactly supported wavelets, using a method tailor-made for the purpose
in \cite{BJKW97}. We will show that all representations obtained from
compactly supported wavelets have a finite-dimensional commutant, and as a
consequence they decompose into a finite direct sum of irreducible
representations. We also display a one-parameter family (with two singular
points) of mutually inequivalent representations of $\mathcal{O}_{2}$ on
$L^{2}\left(  \mathbb{T}\right)  $ for which the corresponding family of
wavelets contains Daubechies's continuous, one-sided differentiable mother
function, $\psi\in L^{2}\left(  \mathbb{R}\right)  $, supported on $\left[
0,3\right]  \subset\mathbb{R}$. In our one-parameter family of wavelets
supported on $\left[  0,3\right]  $, there is actually a left-handed and a
paired right-handed Daubechies wavelet, resulting from a natural symmetry in
the family. In going from one to the other, the one-sided differentiability
property reverses direction.

Let us briefly review how one constructs representations from a
multiresolution wavelet of scale $N$. Many more details may be found in
\cite{BrJo97b}. Excellent accounts of multiresolution wavelet analysis in
general may be found in \cite{Hor95} and \cite{CoRy95}.

Define scaling by $N$ on $L^{2}\left(  \mathbb{R}\right)  $ as the unitary
operator $U$ given by
$\left(  U\xi\right)  \left(  x\right)   =N^{-\frac{1}{2}}\xi\left(
N^{-1}x\right)  $ for $\xi\in L^{2}\left(  \mathbb{R}\right)$
, $x\in\mathbb{R}$,
and translation as the unitary operator $T$ given by
$\left(  T\xi\right)  \left(  x\right)   =\xi\left(  x-1\right)  $.
There is a \emph{father function} or \emph{scaling function}
$\varphi$ which is a vector in $L^{2}\left(
\mathbb{R}\right)  $ such that%
\begin{equation}
\left\{  T^{k}\varphi\right\}  _{k\in\mathbb{Z}}\text{ is an orthonormal set
in }L^{2}\left(  \mathbb{R}\right)  . \label{eqIntro.3}%
\end{equation}
Furthermore, one assumes that there is a sequence $\left(  b_{n}\right)
\in\ell^{2}$ such that%
\begin{equation}
U\varphi=\sum_{n}b_{n}T^{n}\varphi, \label{eqIntro.4}%
\end{equation}
and then necessarily $\sum_{n}\left|  b_{n}\right|  ^{2}=1$. (It seems to be
fairly conventional in wavelet theory to only consider real $b$, but this is
not too important for what follows.) A weaker, so-called ``tight frame'',
property for the vectors in (\ref{eqIntro.3}) will also be considered as a
degenerate case in Section \ref{Casenu=1N=2d=2}. If $\mathcal{V}_{0}$ is the
closed subspace of $L^{2}\left(  \mathbb{R}\right)  $ spanned by $\left\{
T^{k}\varphi\right\}  _{k\in\mathbb{Z}}$, one also assumes%
\begin{equation}
\bigwedge_{n\in\mathbb{Z}}U^{n}\mathcal{V}_{0}  =\left\{  0\right\}
,\qquad
\bigvee_{n\in\mathbb{Z}}U^{n}\mathcal{V}_{0}  =L^{2}\left(  \mathbb{R}
\right)  . \label{eqIntro.5and6}%
\end{equation}

These are all the properties of the father function $\varphi$ that are needed.
One example is the Haar father function $\varphi\left(  x\right)
=\chi_{\left[  0,1\right]  }\left(  x\right)  $.

Define a function $m_{0}$ in $L^{2}\left(  \mathbb{T}\right)  $ by
\begin{equation}
m_{0}\left(  t\right)  =m_{0}\left(  e^{-it}\right)  =\sum_{n}b_{n}e^{-int}.
\label{eqIntro.7}%
\end{equation}
Choose functions $m_{1},\dots,m_{N-1}$ in $L^{2}\left(  \mathbb{T}\right)  $
such that%
\begin{equation}
\sum_{k=0}^{N-1}\overline{m_{i}\left(  t+\frac{2\pi k}{N}\right)  }%
m_{j}\left(  t+\frac{2\pi k}{N}\right)  =\delta_{ij}N \label{eqIntro.8}%
\end{equation}
for almost all $t\in\mathbb{R}$, $i,j=0,1,\dots,N-1$, or, equivalently,
such that the
$N\times N$ matrix%
\begin{equation}
\frac{1}{\sqrt{N}}%
\begin{pmatrix}
m_{0}(z) & m_{0}(\rho z) & \dots &  m_{0}(\rho^{N-1}z)\\
m_{1}(z) & m_{1}(\rho z) & \dots &  m_{1}(\rho^{N-1}z)\\
\vdots & \vdots & \ddots & \vdots\\
m_{N-1}(z) & m_{N-1}(\rho z) & \dots &  m_{N-1}(\rho^{N-1}z)
\end{pmatrix}
, \label{eqIntro.9}%
\end{equation}
where $\rho=e^{\frac{2\pi i}{N}}$, is unitary for almost all $z\in\mathbb{T}$.
(With $m_{0}$ given as above, $m_{1},\dots,m_{N-1}$ may
always be so chosen; see, e.g., \cite{BrJo97b}.)
If we define $\psi_{1},\dots,\psi_{N-1}\in L^{2}\left(  \mathbb{R}\right)  $
by%
\begin{equation}
\sqrt{N}\hat{\psi}_{i}\left(  Nt\right)  =m_{i}\left(  t\right)  \hat{\varphi
}\left(  t\right)  \label{eqIntro.10}%
\end{equation}
for $t\in\mathbb{R}$, $i=1,\dots,N-1$, where $\hat{\;}$ denotes Fourier
transform, unitarity of the above matrix is equivalent to orthonormality in
$L^{2}\left(  \mathbb{R}\right)  $ of the set%
\begin{equation}
\left\{  U^{n}T^{k}\psi_{i}\right\}  _{n,k\in\mathbb{Z};\,
i=1,\dots,N-1}\;. \label{eqIntro.11}%
\end{equation}
The $\psi_{i}$'s are called the mother functions. If $N=2$, there is only one,
of course.

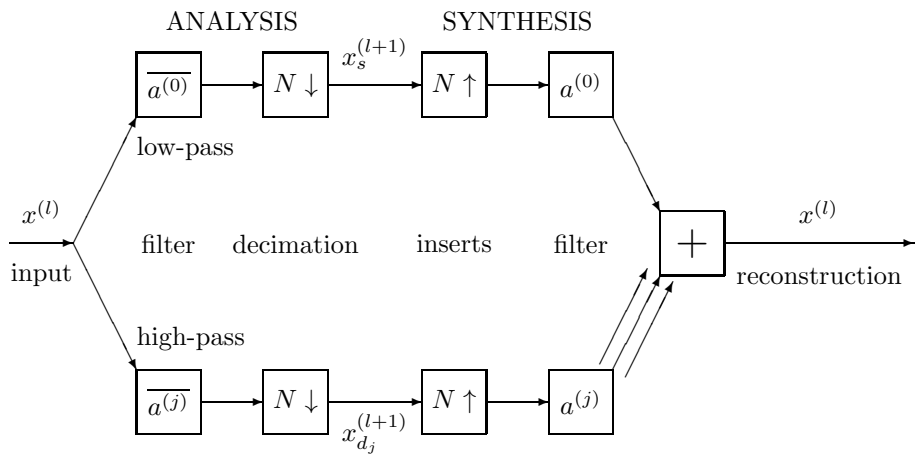
\begin{figure}
\setlength{\unitlength}{2pc}
\begin{picture}(14.5,7.5)(0,-3.5)
\put(0,0){\vector(1,0){1}}
\put(0,0){\makebox(1,1){$x^{(l)}$}}
\put(0,-1){\makebox(1,1){input}}
\put(1,0){\vector(1,2){1}}
\put(1,0){\vector(1,-2){1}}
\put(2,3){\makebox(3,1){ANALYSIS}}
\put(2,2){\framebox(1,1){$\overline{a^{(0)}}$}}
\put(2,1){\makebox(1,1)[l]{low-pass}}
\put(2,-1){\makebox(1,2){filter}}
\put(2,-2){\makebox(1,1)[l]{high-pass}}
\put(2,-3){\framebox(1,1){$\overline{a^{(j)}}$}}
\put(3,2.5){\vector(1,0){1}}
\put(3,-2.5){\vector(1,0){1}}
\put(4,2){\framebox(1,1){$N\downarrow$}}
\put(4,-1){\makebox(1,2){decimation}}
\put(4,-3){\framebox(1,1){$N\downarrow$}}
\put(5,2.5){\vector(1,0){1.5}}
\put(5,2.5){\makebox(1.5,1){$x_{s}^{(l+1)}$}}
\put(5,-2.5){\vector(1,0){1.5}}
\put(5,-3.5){\makebox(1.5,1){$x_{d_{j}}^{(l+1)}$}}
\put(6.5,3){\makebox(3,1){SYNTHESIS}}
\put(6.5,2){\framebox(1,1){$N\uparrow$}}
\put(6.5,-1){\makebox(1,2){inserts}}
\put(6.5,-3){\framebox(1,1){$N\uparrow$}}
\put(7.5,2.5){\vector(1,0){1}}
\put(7.5,-2.5){\vector(1,0){1}}
\put(8.5,2){\framebox(1,1){$a^{(0)}$}}
\put(8.5,-1){\makebox(1,2){filter}}
\put(8.5,-3){\framebox(1,1){$a^{(j)}$}}
\put(9.5,2){\vector(1,-2){0.75}}
\put(9.5,-2){\vector(1,2){0.75}}
\put(9.7,-2.1){\vector(1,2){0.75}}
\put(9.3,-1.9){\vector(1,2){0.75}}
\put(10.25,-0.5){\framebox(1,1){\huge$+$}}
\put(11.25,0){\vector(1,0){3}}
\put(11.25,0){\makebox(3,1){$x^{(l)}$}}
\put(11.25,-1){\makebox(3,1){reconstruction}}
\end{picture}
\caption{Signal subband filtering}
\label{SubbandFiltering}
\end{figure}

Unitarity of (\ref{eqIntro.9}) is also equivalent to saying that the operators
$S_{i}$, defined on $L^{2}\left(  \mathbb{T}\right)  $ by%
\begin{equation}
\left(  S_{i}\xi\right)  \left(  z\right)  =m_{i}\left(  z\right)  \xi\left(
z^{N}\right)  \label{eqIntro.12}%
\end{equation}
for $\xi\in L^{2}\left(  \mathbb{T}\right)  $, $z\in\mathbb{T}$,
$i=0,1,\dots,N-1$, satisfy the relations%
\begin{equation}
S_{j}^{\ast}S_{i}^{{}}  =\delta_{ij}^{{}}\openone,\qquad
\sum_{i=0}^{N-1}S_{i}^{{}}S_{i}^{\ast}  =\openone, \label{eqIntro.13and14}%
\end{equation}
which are exactly the Cuntz relations. There is a one-to-one correspondence
between operator solutions to (\ref{eqIntro.13and14}) and
representations of $\mathcal{O}_{N}$, and since $\mathcal{O}_{N}$ is simple,
these representations are always faithful. The Fourier transform of
$S_{i}^{\ast}$ (the adjoint of (\ref{eqIntro.12})), acting on $\ell^{2}\left(
\mathbb{Z}\right)  $, is the quadrature mirror filter $F_{i}$ in \cite{Wic93}:
$F_{0}$ is low-pass, and $F_{1},\dots,F_{N-1}$ are the corresponding high-pass
filters for the signal reconstitution process. Let
\begin{equation}
m_{i}^{{}}\left(  z\right)  =\sum_{n}a_{n}^{(i)}z_{{}}^{n} \label{eqFoudecomp}%
\end{equation}
be the Fourier decomposition. It follows from (\ref{eqIntro.12}) that for
$x=\left(  x_{k}\right)  _{k\in\mathbb{Z}}\in\ell^{2}$, we have
\begin{equation}
\left(  F_{j}^{\ast}x\right)  _{n}  =\sum_{k\in\mathbb{Z}}a_{n-Nk}^{\left(
j\right)  }x_{k\vphantom{N}}^{{}},\qquad
\left(  F_{j}x\right)  _{n}  =\sum_{k\in\mathbb{Z}}\overline{a_{k-Nn}%
^{\left(  j\right)  }}x_{k\vphantom{N}}^{{}}, \label{eqQuadmf}%
\end{equation}
as operators $\ell^{2}\rightarrow\ell^{2}$. The Cuntz relations in $\ell^{2}%
$-operator form,
\begin{equation}
F_{i}^{{}}F_{j}^{\ast}  =\delta_{ij}^{{}}\openone,\qquad
\sum_{j=0}^{N-1}F_{j}^{\ast}F_{j}^{{}}  =\openone, \label{eqSubband}%
\end{equation}
then summarize subband filtering, which can be written in diagram form as in
Figure \ref{SubbandFiltering}. 
Here ``analysis'' is splitting into subbands and the application of $F_{i}$,
and ``synthesis'' is the application of $F_{i}^{*}$ followed by summing over
the subbands again. The low-pass subband corresponds to $i=0$, and the
high-pass subbands correspond to $i=1,\dots,N-1$.
See \cite{Wic93} and \cite{CoRy95} for details.%

The $\mathcal{O}_{N}$-representations given in (\ref{eqIntro.12}) play a
crucial role in the wavelet analysis in a second related way. A scale-$N$
wavelet in $L^{2}\left(  \mathbb{R}\right)  $ is an orthonormal basis (or a
tight frame) of the form (\ref{eqIntro.11}) as described above. An important
point is that the corresponding $S_{i}$-operators of (\ref{eqIntro.12}), which
constitute the $\mathcal{O}_{N}$-representation, enter directly and explicitly
into a formula for the $L^{2}\left(  \mathbb{R}\right)  $-expansion
coefficients $c_{nki}$ of $\xi=\sum_{n,k,i}c_{nki}\left(  \xi\right)
U^{n}T^{k}\psi_{i}$, $\xi\in L^{2}\left(  \mathbb{R}\right)  $, and we refer
to \cite[eq.~(1.35)]{BrJo97b} for details on that.

\begin{figure}
\mbox{\psfig
{figure=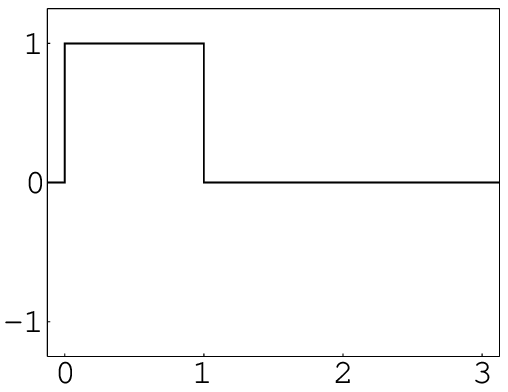,bbllx=0bp,bblly=16bp,bburx=144bp,bbury=128bp,width=120pt}%
\hspace*{40pt}\psfig
{figure=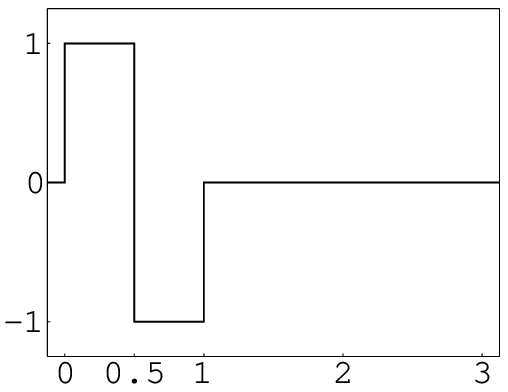,bbllx=0bp,bblly=16bp,bburx=144bp,bbury=128bp,width=120pt}}
\mbox{\makebox[120pt]
{\hspace*{13pt}$\varphi$}\hspace*{40pt}\makebox[120pt]
{\hspace*{13pt}$\psi$}}
\caption{Father and mother functions for the Haar wavelet}
\label{HaarMother}
\end{figure}

We see from (\ref{eqIntro.10}) and (\ref{eqIntro.4}) that the scaled vectors
$U\psi_{i}$ and $U\varphi$ are both finite linear combinations of translates
$\left\{  T^{k}\varphi\right\}  _{k\in\mathbb{Z}}$ if and only if the
functions $m_{i}$ are polynomials, and this is reflected in the fact that the
wavelets $\varphi$, $\psi_{i}$ have compact support if and only if all the
functions $m_{i}\left(  z\right)  $ are polynomials in $z$. (See \cite[Chapter
5]{Dau92}, \cite[Section 3.3]{Hor95}.) In \cite{BrJo96b}, a detailed study was
made of the representations of $\mathcal{O}_{N}$ defined by (\ref{eqIntro.12})
in the case where $m_{i}\left(  z\right)  $ are monomials (or more precisely,
monomials of the form $m_{i}\left(  z\right)  =z^{n_{i}}$; the more general
case where $m_{i}\left(  z\right)  =\lambda_{i}z^{n_{i}}$ with $\lambda_{i}%
\in\mathbb{T\subset C}$ was considered in \cite{DaPi97b}). It is clear from
(\ref{eqIntro.4}) and (\ref{eqIntro.10}) that the other $m_{i}$-functions
coming from wavelets are never monomials, but the Haar wavelet (for $N=2$),
$\varphi\left(  x\right)  =\chi_{\left[  0,1\right]  }\left(  x\right)  $,
is close: one checks from (\ref{eqIntro.4}) and (\ref{eqIntro.7}) that
$m_{0}\left(  z\right)  =\left( 1+z\right) /\sqrt{2}$.
The most general choice of $m_{1}$ is then%
\begin{equation}
m_{1}\left(  z\right)  =zf\left(  z^{2}\right)  \overline{m_{0}\left(
-z\right)  }, \label{eqIntro.17}%
\end{equation}
where $f$ maps $\mathbb{T}$ into $\mathbb{T}$, and one conventional choice is
$f=-1$, i.e.,%
\begin{equation}
m_{1}\left(  z\right)  =\left( 1-z\right) /\sqrt{2}. \label{eqIntro.18}%
\end{equation}
Thus the Haar mother function is given by $\frac{1}{\sqrt{2}}\psi\left(
\frac{x}{2}\right)  =\frac{1}{\sqrt{2}}\left(  \varphi\left(  x\right)
-\varphi\left(  x-1\right)  \right)  $, i.e., the graph of $\psi$ is that
represented in Figure \ref{HaarMother}.%

If $S_{i}$ is defined by (\ref{eqIntro.12}), and one transforms the
representation by $\frac{1}{\sqrt{2}}\left(
\begin{smallmatrix}
1 & -1\\
1 & 1
\end{smallmatrix}
\right)  \in\mathrm{U}\left(  2\right)  $, i.e.,%
\begin{equation}
T_{0}  =\left(  S_{0}+S_{1}\right)  /\sqrt{2},\qquad
T_{1}  =\left(  S_{0}-S_{1}\right)  /\sqrt{2}, \label{eqIntro.19and20}%
\end{equation}
one verifies that the pair $T_{0}$, $T_{1}$ still satisfies the Cuntz
relations, and%
\begin{equation}
T_{0}\xi\left(  z\right)   =\xi\left(  z^{2}\right)  ,\qquad
T_{1}\xi\left(  z\right)   =z\xi\left(  z^{2}\right)  . \label{eqIntro.21and22}%
\end{equation}
This is one of the monomial representations studied in \cite{BrJo96b}, and by
\cite[Proposition 8.1]{BrJo96b}, this representation of $\mathcal{O}_{2}$
decomposes into two inequivalent irreducible subrepresentations on the
subspaces%
\begin{align}
H^{2}\left(  \mathbb{T}\right)   &  =\overline{\operatorname*{span}}^{\left\|
\,\cdot\,\right\|  _{2}}\left\{  z^{n}\mid n\in\mathbb{N}\cup\left\{
0\right\}  \right\} ,\label{eqIntro.23}\\
H^{2}\left(  \mathbb{T}\right)  ^{\perp}  &  =\overline{zH^{2}\left(
\mathbb{T}\right)  }=\overline{\operatorname*{span}}^{\left\|  \,\cdot
\,\right\|  _{2}}\left\{  z^{-n}\mid n\in\mathbb{N}\right\}  ,
\label{eqIntro.24}%
\end{align}
where $\overline{zH^{2}\left(  \mathbb{T}\right)  }$ refers to complex
conjugation, and $\overline{\phantom{\operatorname*{span}}}^{\left\|
\,\cdot\,\right\|  _{2}}$ is $L^{2}\left(  \mathbb{T}\right)  $ closure. Thus
the original Haar wavelet representation is a direct sum of two inequivalent
irreducible subrepresentations. In general, when the functions $m_{i}$ are
polynomials, this simple trick of reducing to monomials is not going to work,
but we will see that it is possible to develop a theory for polynomial
representations which nonetheless has many general analogues with the monomial theory.

\section{\label{Finitely}Finitely correlated states on the Cuntz algebra
$\mathcal{O}_{N}$}

Let us recall a few facts about the Cuntz algebra $\mathcal{O}_{N}$ from
\cite{Cun77}, and the part of the results from \cite{BJKW97} that will be
needed in the sequel.

If $N\in\left\{  2,3,\dots\right\}  $, the \emph{Cuntz algebra} $\mathcal{O}%
_{N}$ is the universal $C^{\ast}$-algebra generated by elements $s_{0}%
,\dots,s_{N-1}$ subject to the relations
\begin{equation}
s_{i}^{\ast}s_{j}^{{}}  =\delta_{ij}\openone,\qquad
\sum_{j\in\mathbb{Z}_{N}}s_{j}^{{}}s_{j}^{\ast}  =\openone.
\label{eqFinitely.1and2}%
\end{equation}

The Cuntz algebra may be
viewed as an interpolation between
the algebra of the canonical anti-commutation
relations (CAR) and
the algebra of the canonical
commutation relations (CCR)%
: The $q$-canonical
commutation relations,
\[
a_{i}^{{}}a_{j}^{*}-qa_{j}^{*}a_{i}^{{}}=\delta_{ij}^{{}}\openone,
\]
$i,j=1,\dots,d$,
reduce to the CCR relations if $q=1$,
the CAR relations if $q=-1$, and
the Cuntz relations (\ref{eqFinitely.1and2}) if $q=0$.
See \cite{JSW94a,BrRoI,BrRoII,
Eva80,
EvKa98} for details on this.

The Cuntz algebra is
a simple separable $C^*$-algebra
not isomorphic to the algebra of
compact operators on a Hilbert
space. Therefore
the space of unitary equivalence
classes of irreducible representations
of $\mathcal{O}_{N}$ cannot be parametrized in
a measurable way \cite{Dix69}. In this
paper we will show that the representations
coming from low-pass filters of genus
$2$ form a (necessarily tiny) one-%
dimensional variety in this enormous space.

There is a canonical action of the group $U\left(  N\right)  $ of unitary
$N\times N$ matrices on $\mathcal{O}_{N}$ given by
\begin{equation}
\tau_{g}\left(  s_{i}\right)  =\sum_{j\in\mathbb{Z}_{N}}\overline{g_{ji}}s_{j}
\label{eqFinitely.3}%
\end{equation}
for $g=\left[  g_{ij}\right]  \in\mathrm{U}\left(  N\right)  $. In particular
the \emph{gauge action} is defined by
$\tau_{z}\left(  s_{i}\right)  =zs_{i}$, $z\in\mathbb{T}\subset\mathbb{C}
\mkern2mu$.
If $\operatorname*{UHF}\nolimits_{N}$ is the fixed point subalgebra under the
gauge action, then $\operatorname*{UHF}\nolimits_{N}$ is the closure of the
linear span of all Wick ordered monomials of the form
$s_{i_{1}}^{{}}\cdots s_{i_{k}}^{{}}s_{j_{k}}^{\ast}\cdots s_{j_{1}}^{\ast}$.
$\operatorname*{UHF}\nolimits_{N}$ is isomorphic to the $\operatorname*{UHF}%
$-algebra of Glimm type $N^{\infty}$,
\begin{equation}
\operatorname*{UHF}\nolimits_{N}\cong M_{N^{\infty}}=\bigotimes_{1}^{\infty
}M_{N}, \label{eqFinitely.6}%
\end{equation}
in such a way that the isomorphism carries the aforementioned
Wick ordered monomial,
$s_{i_{1}}^{{}}\cdots s_{i_{k}}^{{}}s_{j_{k}}^{\ast}\cdots s_{j_{1}}^{\ast}$,
into the matrix element%
\begin{equation}
e_{i_{1}j_{1}}\otimes e_{i_{2}j_{2}}\otimes\dots\otimes e_{i_{k}j_{k}}%
\otimes\openone\otimes\openone\otimes\cdots. \label{eqFinitely.7}%
\end{equation}
The restriction of $\tau_{g}$ to $\operatorname*{UHF}\nolimits_{N}$ is then
carried into the action%
\begin{equation}
\operatorname{Ad}\left(  g\right)  \otimes\operatorname{Ad}\left(  g\right)
\otimes\cdots\label{eqFinitely.8}%
\end{equation}
on $\bigotimes_{1}^{\infty}M_{N}$. We define the canonical endomorphism
$\lambda$ on $\operatorname*{UHF}\nolimits_{N}$ (or on $\mathcal{O}_{N}$) by
\begin{equation}
\lambda\left(  x\right)  =\sum_{j\in\mathbb{Z}_{N}}s_{j}^{{}}xs_{j}^{\ast}
\label{eqFinitely.9}%
\end{equation}
and the isomorphism carries $\lambda$ over into the one-sided shift
\begin{equation}
x_{1}\otimes x_{2}\otimes x_{3}\otimes\cdots\longrightarrow\openone\otimes
x_{1}\otimes x_{2}\otimes\cdots\label{eqFinitely.10}%
\end{equation}
on $\bigotimes_{1}^{\infty}M_{N}$.

If $s_{i}\mapsto S_{i}\in\mathcal{B}\left(  \mathcal{H}\right)  $ is a
representation of the Cuntz relations on a Hilbert space $\mathcal{H}$, we
will say (by abuse of terminology) that the representation is finitely
correlated if there exists a finite-dimensional subspace $\mathcal{K}%
\subset\mathcal{H}$ with the two properties
\begin{gather}
\begin{minipage}[t]{\displayboxwidth}\raggedright$\displaystyle S_{i}^{\ast
}\mathcal{K}\subset\mathcal{K}$\quad for $\displaystyle i\in\mathbb{Z}%
_{N}$, \end{minipage} \label{eqFinitely.11}\\
\begin{minipage}[t]{\displayboxwidth}\raggedright$\displaystyle\mathcal
{K}$ is cyclic for the representation $\displaystyle s_{i}\longmapsto
S_{i}$. \end{minipage} \label{eqFinitely.12}%
\end{gather}
The presence of such a finite-dimensional subspace $\mathcal{K}$ is a special
property of each of the representations under discussion, and therefore of the
states of $\mathcal{O}_{N}$ which correspond to the representations. These
states were studied in \cite{BJKW97} with a view to the present applications.

If $P\colon\mathcal{H}\rightarrow\mathcal{K}$ is the orthogonal projection
onto $\mathcal{K}$, then (\ref{eqFinitely.11}) can be formulated as
\begin{equation}
V_{i}\equiv PS_{i}=PS_{i}P. \label{eqFinitely.13}%
\end{equation}
If we view $V_{i}$ as operators in $\mathcal{B}\left(  \mathcal{K}\right)  $,
we have
\begin{equation}
\sum_{i\in\mathbb{Z}_{N}}V_{i}^{{}}V_{i}^{\ast}=\openone,
\label{eqFinitely.14}%
\end{equation}
and conversely, if $V_{i}$ are operators in $\mathcal{B}\left(  \mathcal{K}%
\right)  $ satisfying (\ref{eqFinitely.14}), they determine a representation
$s_{i}\mapsto S_{i}$ of the Cuntz relations such that (\ref{eqFinitely.13}) is
valid, and this representation is unique up to unitary equivalence if we
require $\mathcal{K}$ to be cyclic \cite[Theorem 5.1]{BJKW97}.

If $\mathcal{K}_{1}$ is another Hilbert space and $W_{0},\dots,W_{N-1}$ are
operators on $\mathcal{K}_{1}$ satisfying%
\begin{equation}
\sum_{i\in\mathbb{Z}_{N}}W_{i}^{{}}W_{i}^{\ast}=\openone,
\label{eqFinitely.15}%
\end{equation}
and $s_{i}\mapsto T_{i}$ is the associated representation of $\mathcal{O}_{N}%
$, then there is an isometric linear isomorphism between intertwiners
$U\colon\mathcal{H}_{V}\rightarrow\mathcal{H}_{W}$, i.e., operators
satisfying
\begin{equation}
US_{i}=T_{i}U, \label{eqFinitely.16}%
\end{equation}
and operators $V\in\mathcal{B}\left(  \mathcal{K},\mathcal{K}_{1}\right)  $
such that
\begin{equation}
\br\left(  V\right)  \equiv\sum_{i\in\mathbb{Z}_{N}}W_{i}^{{}}VV_{i}^{\ast}=V.
\label{eqFinitely.17}%
\end{equation}
This linear isomorphism is given by
\begin{equation}
U\longmapsto V=P_{1}VP, \label{eqFinitely.18}%
\end{equation}
where $P_{1}\colon\mathcal{H}_{W}\rightarrow\mathcal{K}_{1}$ is the orthogonal
projection onto $\mathcal{K}_{1}$. All these results do not depend on
$\mathcal{K}$ and $\mathcal{K}_{1}$ being finite-dimensional, and they are
given in \cite[Theorem 5.1]{BJKW97}.

An important special case is $\mathcal{K}_{1}=\mathcal{K}$ and $W_{i}=V_{i}$.
Then $\br$ is a completely positive unital map, and the linear isomorphism
(\ref{eqFinitely.18}) is an order isomorphism between the fixed point set of
$\br$ (which is not necessarily an algebra) and the commutant $\left\{
S_{i}^{{}},S_{i}^{\ast}\mid i\in\mathbb{Z}_{N}\right\}  ^{\prime}$. In
particular, we have the following principle.%
\begin{equation}
\begin{minipage}[t]{\displayboxwidth}\raggedright
The representation $\displaystyle s_{i}\longmapsto S_{i}%
$ is irreducible if and only if $\displaystyle\br$ is ergodic: $\displaystyle
\left\{ A\in\mathcal{B}\left( \mathcal{K}\right) \mid\br\left( A\right
) =A\right\} =\mathbb{C}\mkern2mu\openone$. \end{minipage}
\label{eqFinitely.19}%
\end{equation}
The rest of the discussion in this section can only be partially extended to
the case when $\mathcal{K}$ is infinite-dimensional (see \cite[Section
6]{BJKW97} for details). Define $\bs=\br$ in the case when $\mathcal{K}%
_{1}=\mathcal{K}$ and $W_{i}=V_{i}$ in (\ref{eqFinitely.17}). If $\bs$ is
ergodic, then $\mathcal{B}\left(  \mathcal{K}\right)  $ has a unique $\bs
$-invariant state $\varphi$. This state need not be faithful (see the example
after the proof of Lemma 3.4 in \cite{BJKW97}). If $E$ is the support
projection of $\varphi$, then $S_{i}^{\ast}E\mathcal{K}\subset E\mathcal{K}$
for all $i\in\mathbb{Z}_{N}$ (see \cite[Lemma 6.1]{BJKW97}). In that case,
replace $P$ by $E$, $V_{i}$ by $EV_{i}$, $\bs$ by the $\bs$ defined by the new
$V_{i}$'s on $E\mathcal{K}$, and then define a state $\psi$ on $\mathcal{O}%
_{N}$ by%
\begin{equation}
\psi\left(  S_{I}^{{}}S_{J}^{\ast}\right)  =\varphi\left(  ES_{I}^{{}}%
S_{J}^{\ast}E\right)  . \label{eqFinitely.20}%
\end{equation}
It was proved in \cite[Theorem 6.3]{BJKW97} that the following three subsets
of the circle group $\mathbb{T}$ are equal:%
\begin{gather}
\begin{minipage}[t]{\displayboxwidth}\raggedright$\displaystyle\left
\{ t\in\mathbb{T}\mid\psi\circ\tau_{t} =\psi\right\} $, where $\displaystyle
\tau$ is the gauge action; \end{minipage} \label{eqFinitely.21}\\
\begin{minipage}[t]{\displayboxwidth}\raggedright$\displaystyle\left
\{ t\in\mathbb{T}\mid\psi\circ\tau_{t}\text{ is quasi-equivalent to }%
\psi\right\} $; \end{minipage} \label{eqFinitely.22}\\
\begin{minipage}[t]{\displayboxwidth}\raggedright$\displaystyle\operatorname
{PSp}\left( \bs\right) \cap\mathbb{T}$, where $\displaystyle\operatorname
{PSp}\left( \bs\right) $ is the set of eigenvalues of $\displaystyle
\bs$. \end{minipage} \label{eqFinitely.23}%
\end{gather}
(Of course, in the present setting, where $E\mathcal{K}$ is
finite-dimensional, $\operatorname*{PSp}\left(  \bs\right)
=\operatorname*{Sp}\left(  \bs\right)  $.) Furthermore, this subset is a
finite subgroup of $\mathbb{T}$. If $k$ is the order of this subgroup, the
restriction of the representation to $\operatorname*{UHF}\nolimits_{N}$
decomposes into $k$ mutually disjoint irreducible representations, and these
are mapped cyclically one into another by the one-sided shift $\lambda$. More
specifically, one has
$\operatorname{PSp}\left(  \bs\right)  \cap\mathbb{T}=\operatorname{PSp}\left(
\lambda\right)  \cap\mathbb{T}$,
and, if $t_{k}=e^{\frac{2\pi i}{k}}$, there exists a unitary $U$ on
$\mathcal{H}$, unique up to a scalar, implementing $\tau_{t_{k}}$, and such
that $U^{k}=\openone$. The operator $U$ is the unique (up to a scalar)
eigen-element such that
$\lambda\left(  U\right)  =\bar{t}_{k}U$.
If%
\begin{equation}
U=\sum_{n\in\mathbb{Z}_{k}}t_{k}^{n}E_{n}^{{}} \label{eqFinitely.26}%
\end{equation}
is the spectral decomposition of $U$, then the spectral projections $E_{n}$
project into mutually disjoint irreducible subspaces invariant for the
representation restricted to $\operatorname*{UHF}\nolimits_{N}$, and
$\lambda\left(  E_{n}\right)  =E_{n+1}$, with $\lambda$ extended to
$\mathcal{B}\left(  \mathcal{H}\right)  $ by the formula
$\lambda\left(  \,\cdot\,\right)  =\sum_{i\in\mathbb{Z}_{N}}S_{i}^{{}}
\,\cdot\,S_{i}^{\ast}$.

\section{\label{Polynomial}Polynomial representations}

{}From the relation (\ref{eqIntro.12}) it follows that%
\begin{equation}
\left(  S_{i}^{\ast}\xi\right)  \left(  z\right)  =\frac{1}{N}\sum_{w^{N}%
=z}\overline{m_{i}\left(  w\right)  }\xi\left(  w\right)  ,
\label{eqPolynomial.1}%
\end{equation}
where the sum ranges over all $N$'th roots $w$ of $z$ \cite[eq.~%
(1.17)]{BrJo97b}. Recall that the Fourier series version of
(\ref{eqPolynomial.1}) on $\ell^{2}\left(  \mathbb{Z}\right)  $ is the filter
operator $F_{i}$ of (\ref{eqQuadmf}). In order to incorporate the monomial
results obtained in \cite{BrJo96b}, and also to make the present results
applicable to wavelets in dimension $\nu>1$, let us extend the definitions of
the representations somewhat. We replace $L^{2}\left(  \mathbb{T}\right)  $
with $L^{2}\left(  \mathbb{T}^{\nu}\right)  $ and fix a matrix $\mathbf{N}$
with integer coefficients such that $\left|  \det\left(  \mathbf{N}\right)
\right|  =N\in\left\{  2,3,\dots\right\}  $. If $z=\left(  z_{1},\dots,z_{\nu
}\right)  \in\mathbb{T}^{\nu}$ define%
\begin{equation}
z^{\mathbf{N}}=\left(  z_{1}^{n_{11}}\cdots z_{\nu}^{n_{\nu1}},\dots
,z_{1}^{n_{1\nu}}\cdots z_{\nu}^{n_{\nu\nu}}\right)  \in\mathbb{T}^{\nu}
\label{eqPolynomial.2}%
\end{equation}
if $\mathbf{N}=\left[  n_{ij}\right]  _{i,j=1}^{\nu}$. (Note that this
definition of $z^{\mathbf{N}}$ is different from the one after (1.8) in
\cite{BrJo96b}. The present convention implies that relations like $\left(
z^{\mathbf{N}}\right)  ^{\mathbf{M}}=z^{\mathbf{NM}}$ and $\left(
z^{\mathbf{N}}\right)  ^{n}=z^{\mathbf{N}n}$ are valid, where $z^{n}$ is
defined
as in connection with (\ref{eqPolynomial.6}) below.
The present map $z\mapsto z^{\mathbf{N}}$ is the transpose of the
map $x\mapsto \mathbf{N}x$ on $\mathbb{R}^{\nu}$ passed to the quotient
$\mathbb{T}^{\nu}=\mathbb{R}^{\nu}\diagup 2\pi\mathbb{Z}^{\nu}$.)
The map $z\mapsto z^{\mathbf{N}}$ is
$N$-to-$1$. Let $\sigma_{0},\dots,\sigma_{N-1}$ denote sections of this map,
i.e., each $\sigma_{i}\colon\mathbb{T}^{\nu}\rightarrow\mathbb{T}^{\nu}$ is
injective, $\mu\left(  \sigma_{i}\left(  \mathbb{T}^{\nu}\right)  \cap
\sigma_{j}\left(  \mathbb{T}^{\nu}\right)  \right)  =0$ if $i\neq j$, where
$\mu$ is normalized Haar measure on $\mathbb{T}^{\nu}$, and $\mu\left(
\sigma_{i}\left(  Y\right)  \right)  =\frac{1}{N}\mu\left(  Y\right)  $ for
all Borel sets $Y\subset\mathbb{T}^{\nu}$. Thus $\bigcup_{i\in\mathbb{Z}_{N}%
}\sigma_{i}\left(  \mathbb{T}^{\nu}\right)  =\mathbb{T}^{\nu}$ up to sets of
measure zero. The unitarity condition (\ref{eqIntro.9}) then says that the
$N\times N$ matrix%
\begin{equation}
\frac{1}{\sqrt{N}}%
\begin{pmatrix}
m_{0}(\sigma_{0}(z)) & m_{0}(\sigma_{1}(z)) & \dots &  m_{0}(\sigma
_{N-1}(z))\\
m_{1}(\sigma_{0}(z)) & m_{1}(\sigma_{1}(z)) & \dots &  m_{1}(\sigma
_{N-1}(z))\\
\vdots & \vdots & \ddots & \vdots\\
m_{N-1}(\sigma_{0}(z)) & m_{N-1}(\sigma_{1}(z)) & \dots &  m_{N-1}%
(\sigma_{N-1}(z))
\end{pmatrix}
\label{eqPolynomial.3}%
\end{equation}
is unitary for almost all $z\in\mathbb{T}^{\nu}$. The representation
(\ref{eqIntro.12}), (\ref{eqPolynomial.1}) of $\mathcal{O}_{N}$ now takes the
form%
\begin{align}
\left(  S_{i}\xi\right)  \left(  z\right)   &  =m_{i}\left(  z\right)
\xi\left(  z^{\mathbf{N}}\right)  ,\label{eqPolynomial.4}\\%
\intertext{and then}%
\left(  S_{i}^{\ast}\xi\right)  \left(  z\right)   &  =\frac{1}{N}%
\sum_{w^{\mathbf{N}}=z}\overline{m_{i}\left(  w\right)  }\xi\left(  w\right)
. \label{eqPolynomial.5}%
\end{align}

Now, assume in addition to unitarity of (\ref{eqPolynomial.3}) that
$m_{0},\dots,m_{N-1}$ all are polynomials, so that there exists a fixed finite
subset $D\subset\mathbb{Z}^{\nu}$ such that%
\begin{equation}
m_{j}\left(  z\right)  =\sum_{n\in D}a_{n}^{\left(  j\right)  }z^{n}.
\label{eqPolynomial.6}%
\end{equation}
Here we have used the notation
$z^{n}=\left(  z_{1},\dots,z_{\nu}\right)  ^{\left(  n_{1},\dots,n_{\nu
}\right)  }=z_{1}^{n_{1}}z_{2}^{n_{2}}\cdots z_{\nu}^{n_{\nu}}$,
and $a_{n}^{\left(  j\right)  }\in\mathbb{C}%
\mkern2mu%
$. Let $e_{n}$, $n\in\mathbb{Z}^{\nu}$, denote the usual Fourier basis for
$L^{2}\left(  \mathbb{Z}^{\nu}\right)  $, i.e.,
$e_{n}\left(  z\right)  =z^{n}$.
It follows from (\ref{eqPolynomial.4}) that
\begin{equation}
S_{j}e_{n}=\sum_{k\in D}a_{k}^{\left(  j\right)  }e_{k+\mathbf{N}n}^{{}}.
\label{eqPolynomial.9}%
\end{equation}
If in general we define $a_{k}^{\left(  j\right)  }=0$ when $k\notin D$, it
follows from (\ref{eqPolynomial.9}) or (\ref{eqPolynomial.5}) that%
\begin{equation}
S_{j}^{\ast}e_{n}^{{}}=\sum_{m\in\mathbb{Z}^{\nu}}\overline{a_{n-\mathbf{N}%
m}^{\left(  j\right)  }}e_{m\vphantom{\mathbf{N}}}^{{}}=\sum_{p\in
D\colon p=n\operatorname{mod}\mathbf{N}}\overline{a_{p\vphantom{\mathbf{N}}%
}^{\left(  j\right)  }}e_{\mathbf{N}^{-1}\left(  n-p\right)  }^{{}}.
\label{eqPolynomial.10}%
\end{equation}
Thus both $S_{j}^{{}}$ and $S_{j}^{\ast}$ map trigonometric polynomials into
trigonometric polynomials in this case. If the matrix $\mathbf{N}^{-1}$
defines a contractive map $\mathbb{R}^{\nu}\rightarrow\mathbb{R}^{\nu}$ in
some norm, one can say more. The following proposition is an analogue of Lemma
3.8 in \cite{BrJo96b} in the present setting.

\begin{proposition}
\label{ProPolynomial.1}Assume that all the \textup{(}complex\textup{)}
eigenvalues of $\mathbf{N}$ have modulus greater than $1$. It follows that
there is a finite subset $H\subset\mathbb{Z}^{\nu}$ with the property that for
any $n\in\mathbb{Z}^{\nu}$ there exists an $M\in\mathbb{N}$ such that%
\begin{equation}
S_{I}^{\ast}e_{n}^{{}}\in\widehat{\ell^{2}\left(  H\right)  }\equiv
\operatorname*{span}\left\{  e_{m}\mid m\in H\right\}  \label{eqPolynomial.11}%
\end{equation}
for all multi-indices $I$ with $\left|  I\right|  \geq M$.
\end{proposition}

\begin{proof}
Let us give two proofs of this statement, both based on a study of the maps
$\sigma_{p}\colon\mathbb{R}^{\nu}\rightarrow\mathbb{R}^{\nu}$ defined for
$p\in D$ by%
\begin{equation}
\sigma_{p}\left(  x\right)  =\mathbf{N}^{-1}\left(  x-p\right)
\label{eqPolynomial.12}%
\end{equation}
for $x\in\mathbb{R}^{\nu}$. By considering a Jordan form of $\mathbf{N}$, as
in the proof of Lemma 3.8 in \cite{BrJo96b}, the condition $\left|
\lambda_{i}\right|  >1$ on the eigenvalues of $\mathbf{N}$ means that there
exists a norm on $\mathbb{C}^{\nu}$ such that
$\left\|  \mathbf{N}^{-1}\right\|  <1$
in the associated norm on $\mathcal{B}\left(  \mathbb{C}^{\nu}\right)  $. If
$d=\max\left\{  \left\|  p\right\|  \mid p\in D\right\}  $,
it follows from (\ref{eqPolynomial.12}) that
$\left\|  \sigma_{p}\left(  x\right)  \right\|  \leq\left\|  \mathbf{N}
^{-1}\right\|  \left(  \left\|  x\right\|  +d\right)  $
for $p\in D$, and by iteration,%
\begin{multline}
\left\| \sigma_{p_{1}}\sigma_{p_{2}}\cdots\sigma_{p_{n}}
\left( x\right) \right\| \leq\left\| \mathbf{N}^{-1}\right\| ^{n}
\left\| x\right\| +\sum_{k=1}^{n}\left\| \mathbf{N}^{-1}\right\| ^{k}
d\label{eqPolynomial.16} \\
=\left\| \mathbf{N}^{-1}\right\| ^{n}\left\| x\right\| +\left
\| \mathbf{N}^{-1}\right\| \frac{1-\left\| \mathbf{N}^{-1}\right\| ^{n}
}{1-\left\| \mathbf{N}^{-1}\right\| }d
\leq\left\| \mathbf{N}^{-1}
\right\| ^{n}\left\| x\right\| +\frac{\left\| \mathbf{N}^{-1}\right
\| }{1-\left\| \mathbf{N}^{-1}\right\| }d
\end{multline}
for $p_{1},\dots,p_{n}\in D$, $n\in\mathbb{N}$. Now, using
(\ref{eqPolynomial.10}) in the form%
\begin{equation}
S_{j}^{\ast}e_{n}^{{}}=\sum_{p\in D\colon p=n\operatorname{mod}%
\mathbf{N}}\overline{a_{p\vphantom{\sigma_{p}\left( n\right) }}^{\left(
j\right)  }}e_{\sigma_{p}\left(  n\right)  }^{{}}, \label{eqPolynomial.17}%
\end{equation}
one deduces from (\ref{eqPolynomial.16}) that%
\begin{multline}
S_{I}^{\ast}\widehat{\ell^{2}}\left(  \left\{  m\in\mathbb{Z}^{\nu}%
\mid\left\|  m\right\|  \leq R\right\}  \right) \label{eqPolynomial.18}\\
\subset\widehat{\ell^{2}}\left(  \left\{  m\in\mathbb{Z}^{\nu}\biggm|\left\|
m\right\|  \leq\left\|  \mathbf{N}^{-1}\right\|  ^{\left|  I\right|  }%
R+\left(  \left\|  \mathbf{N}^{-1}\right\|  /\left(  1-\left\|  \mathbf{N}%
^{-1}\right\|  \right)  \right)  d\right\}  \right)  .
\end{multline}
Thus Proposition \ref{ProPolynomial.1} follows with
\begin{equation}
H=\left\{  n\in\mathbb{Z}^{\nu}\biggm|\left\|  n\right\|  \leq
\left(  \left\|  \mathbf{N}^{-1}\right\|  /\left(  1-\left\|  \mathbf{N}%
^{-1}\right\|  \right)  \right)  d\right\}  .
\settowidth{\qedskip}{$\displaystyle
H=\left\{  n\in\mathbb{Z}^{\nu}\mid\left\|  n\right\|  \leq
\left(  \left\|  \mathbf{N}^{-1}\right\|  /\left(  1-\left\|  \mathbf{N}%
^{-1}\right\|  \right)  \right)  d\right\}
.$}
\addtolength{\qedskip}{-\textwidth}
\rlap{\hbox to-0.5\qedskip{\hfil\qed}}
\label{eqPolynomial.19}%
\end{equation}%
\renewcommand{\qed}{}%
\end{proof}

\begin{remark}
\label{RemPolynomial.2}The other method of proving Proposition
\textup{\ref{ProPolynomial.1}} is a small variation which gives an optimal
choice of $H$ given only $D$. By a theorem of Bandt
\textup{\cite{Ban91,Ban97,DDL95,Str94}} cited in \textup{\cite[(3.11)--(3.12)]%
{BrJo96b}} there is a unique compact subset $X\subset\mathbb{R}^{\nu}$ such
that $X$ is a fixed point for the map $Y\mapsto\bigcup_{p\in D}\sigma
_{p}\left(  Y\right)  $, i.e.,%
\begin{equation}
X=\bigcup_{p\in D}\sigma_{p}\left(  X\right)  , \label{eqPolynomial.20}%
\end{equation}
and we may take%
\begin{equation}
H=X\cap\mathbb{Z}^{\nu}=H\left(  D\right)  . \label{eqPolynomial.21}%
\end{equation}
In some examples in Section \textup{\ref{nu=1}}, the finite subset
$H\subset\mathbb{Z}^{\nu}$ will be computed explicitly.
If the representation of $\mathcal{O}_{N}$ is irreducible, an
application of \cite[Lemma 6.1]{BJKW97}
further shows that
the finite-dimensional subspace
$\mathcal{K}\left(  H\right) $ from \textup{(\ref{eqPolynomial.21})}
contains a
unique minimal subspace $\mathcal{M}\neq 0$ with
the invariance property
$S_{i}^{*}\mathcal{M}\subset \mathcal{M}$.
\end{remark}

The following corollary is the main tool in analyzing polynomial representation.

\begin{corollary}
\label{CorPolynomial.3}Consider the polynomial representation of
$\mathcal{O}_{N}$ defined by \textup{(\ref{eqPolynomial.9})} and
\textup{(\ref{eqPolynomial.10})}, and let $H$ be a minimal finite subset of
$\mathbb{Z}^{\nu}$ satisfying the properties in Proposition
\textup{\ref{ProPolynomial.1}}. It follows that%
\begin{equation}
\mathcal{K}=\widehat{\ell^{2}\left(  H\right)  } \label{eqPolynomial.22}%
\end{equation}
is cyclic for the representation, and thus the representation is finitely
correlated. Defining $V_{j}^{\ast}\in\mathcal{B}\left(  \mathcal{K}\right)  $
by%
\begin{equation}
V_{j}^{\ast}e_{n}^{{}}=\sum_{m\in H}\overline{a_{n-\mathbf{N}m}^{\left(
j\right)  }}e_{m\vphantom{\mathbf{N}}}^{{}} =\sum
_{\substack{p\in D\colon p=n\operatorname{mod}\mathbf{N}\\\sigma_{p}\left(
n\right)  \in H}}\overline{a_{p \vphantom{\sigma_{p}\left( n\right) }
}^{\left(  j\right)  }}e_{\sigma_{p}\left(  n\right)  }^{{}}
\label{eqPolynomial.23}%
\end{equation}
for $n\in H$, the commutant of the representation is isometrically order
isomorphic to%
\begin{equation}
\mathcal{B}\left(  \mathcal{K}\right)  ^{\bs}=\left\{  A\in\mathcal{B}\left(
\mathcal{K}\right)  \biggm|\bs\left(  A\right)  \equiv\sum_{k\in\mathbb{Z}%
_{N}}V_{k}^{{}}AV_{k}^{\ast}=A\right\}  . \label{eqPolynomial.24}%
\end{equation}
In particular the representation is irreducible if and only if $\mathcal{B}%
\left(  \mathcal{K}\right)  ^{\bs}=\mathbb{C}%
\mkern2mu%
\openone$. In this case, the peripheral spectrum of $\bs$ is always a finite
(necessarily cyclic) subgroup of $\mathbb{T}$, and if $k$ is the order of this
subgroup, the restriction of the representation to $\operatorname*{UHF}%
\nolimits_{N}$ decomposes into the direct sum of $k$ mutually disjoint
irreducible representations.

In general the intertwiner space between two representations of this type is
given by \textup{(\ref{eqFinitely.16})--(\ref{eqFinitely.17})}.
\end{corollary}

\begin{proof}
The identity%
\begin{equation}
\openone=\sum_{I\colon\left|  I\right|  =M}S_{I}^{{}}S_{I}^{\ast},
\label{eqPolynomial.25}%
\end{equation}
in conjunction with Proposition \ref{ProPolynomial.1}, implies that all
monomials $e_{n}$, $n\in\mathbb{Z}^{\nu}$, are contained in the cyclic
subspace generated by $\mathcal{K}$, and hence this space is dense in
$L^{2}\left(  \mathbb{T}^{\nu}\right)  $. Indeed, for every $n\in
\mathbb{Z}^{\nu}$, there is, by Proposition \ref{ProPolynomial.1}, an
$M\in\mathbb{N}$ such that $S_{I}^{\ast}e_{n}^{{}}\in\mathcal{K}$ for all $I$
such that $\left|  I\right|  \geq M$. Therefore $S_{I}^{{}}S_{I}^{\ast}%
e_{n}^{{}}\in S_{I}^{{}}\left(  \mathcal{K}\right)  $. An application of
(\ref{eqPolynomial.25}) to $e_{n}$ then yields the desired cyclicity. This
cyclicity is the second of the two properties of the subspace $\mathcal{K}$ in
the discussion of Section \ref{Finitely}, i.e., (\ref{eqFinitely.12}). The
rest (and some more details) follows from the discussion in Section
\ref{Finitely}.
\end{proof}

\section{\label{nu=1}Classification of some polynomial representations}

If $D$ is a given finite subset of $\mathbb{Z}^{\nu}$, the set of all
polynomials $m_{j}$ given by (\ref{eqPolynomial.6}), and satisfying the
unitarity condition (\ref{eqPolynomial.3}) and the normalization%
\begin{equation}
m_{0}\left(  1\right)  =\sqrt{N} \label{eqnu=1.1}%
\end{equation}
(which is necessary for the convergence of the Mallat expansion; see
\cite{Mal89} or \cite[eq.~(1.37)]{BrJo97b}), forms a compact algebraic
variety $\mathcal{M}_{D}$, and it is given as the solution variety of a set of
quadratic equations in the coefficients $a_{n}^{\left(  j\right)  }$ and
$\overline{a_{n}^{\left(  j\right)  }}$ with $n\in D$. For each point on this
variety $\mathcal{M}_{D}$, the corresponding representation of $\mathcal{O}%
_{N}$ can in principle be computed from Corollary \ref{CorPolynomial.3}. Even
the characterization of $\mathcal{M}_{D}$ is a formidable task in general, but
it has been done in the case $\nu=1$ and $N=2$ in
\cite{Wel93,ReWe92,McWe94,HeWe94,HRW92} and \cite{Pol89} (see also
\cite{Pol90,Pol92}). In this section, we will compute the representation
theory of $\mathcal{O}_{2}$ for each of the points of some of these varieties.
We do not know if our results indicate how the generic behaviour of this
representation theory will be, but in the examples the representations
generically are irreducible and mutually disjoint, with exceptional behaviour
on a sub-variety of lower dimension.

\subsection{\label{Casenu=1}The case with dimension $\nu=1$}

In this case, $\mathbf{N}=N\in\left\{  2,3,4,\dots\right\}  $. If
$m_{j}\left(  z\right)  =\sum_{n\in D}a_{n}^{\left(  j\right)  }z_{{}}^{n}$,
where $m=m_{i}$ for some $i$, unitarity of (\ref{eqIntro.9}) implies%
\begin{equation}
\sum_{k\in\mathbb{Z}_{N}}\left|  m_{j}\left(  \rho^{k}z\right)  \right|
^{2}=N, \label{eqnu=1.3}%
\end{equation}
which is equivalent to the conditions%
\begin{equation}
\sum_{n}a_{n}^{\left(  j\right)  }\overline{\displaystyle a_{n}^{\left(
j\right)  }}  =1\text{\quad and\quad }
\sum_{n}a_{n\vphantom{N}}^{\left(  j\right)  }\overline{a_{n-mN}^{\left(
j\right)  }}  =0 \label{eqnu=1.4and5}%
\end{equation}
for $m=1,2,\dots$. Analogously, orthogonality of the rows in (\ref{eqIntro.9})
leads to%
\begin{equation}
\sum_{n}a_{n\vphantom{N}}^{\left(  i\right)  }\overline{a_{n-mN}^{\left(
j\right)  }}=0 \label{eqnu=1.6}%
\end{equation}
for all $i\neq j$ and \emph{all} $m\in\mathbb{Z}$. Finally, the normalization
(\ref{eqnu=1.1}) leads to%
\begin{equation}
\sum_{n}a_{n}^{\left(  0\right)  }=\sqrt{N}. \label{eqnu=1.7}%
\end{equation}

The relations (\ref{eqnu=1.4and5})--(\ref{eqnu=1.7}), together with $a_{n}%
^{\left(  i\right)  }=0$ for $n\notin D$, determine the algebraic variety
$\mathcal{M}_{D}$. Let us now restrict to $N=2$, and to the case where the
$a_{n}^{\left(  j\right)  }$'s are real (this latter assumption, reality,
seems conventional in wavelet theory). Then by (\ref{eqIntro.17}),%
\begin{equation}
m_{1}\left(  z\right)  =zf\left(  z^{2}\right)  \overline{m_{0}\left(
-z\right)  }, \label{eqnu=1.8}%
\end{equation}
where $f$ is a monomial. By translating the father, and mother, functions by
multiples of $T$ (integral translations), we may assume that $D$ has the form%
\begin{equation}
D=\left\{  0,1,\dots,2d-1\right\}  , \label{eqnu=1.9}%
\end{equation}
where $d\in\mathbb{N}$, and with $f\left(  w\right)  =-w^{d-1}$, we have%
\begin{align}
m_{0}\left(  z\right)   &  =\sum_{k=0}^{2d-1}a_{k}z^{k},\label{eqnu=1.10}\\
m_{1}\left(  z\right)   &  =\sum_{k=0}^{2d-1}\left(  -1\right)  ^{k+1}%
a_{k}z^{2d-1-k}=\sum_{k=0}^{2d-1}\left(  -1\right)  ^{k}a_{2d-1-k}z^{k}.
\label{eqnu=1.11}%
\end{align}

The conditions (\ref{eqnu=1.4and5})--(\ref{eqnu=1.7}) then become%
\begin{equation}
\sum_{k=0}^{2d-1}a_{k}^{2}  =1,\text{\quad and\quad}
\sum_{k=0}^{2\left(  d-m\right)  -1}a_{k}a_{k+2m}  =0, \label{eqnu=1.12and13}%
\end{equation}
for $m=1,\dots,d-1$ (no condition if $d=1$), and%
\begin{equation}
\sum_{k=0}^{2d-1}a_{k}=\sqrt{2}. \label{eqnu=1.14}%
\end{equation}
(The condition (\ref{eqnu=1.6}) is already taken care of in (\ref{eqnu=1.11}).)

In this case the maps $\sigma_{p}$ in (\ref{eqPolynomial.12}) have the form%
\begin{equation}
\sigma_{p}\left(  x\right)  =\frac{x-p}{2} \label{eqnu=1.15}%
\end{equation}
for $p=0,1,\dots,2d-1$, and thus the solution $X$ to the equation
(\ref{eqPolynomial.20}) is the interval%
\begin{equation}
X  =\left[  -2d+1,0\right]  ,\label{eqnu=1.16}%
\end{equation}
and hence
by (\ref{eqPolynomial.21})
\begin{align}
H  &  =\left\{  -2d+1,-2d+2,\dots,0\right\} ,\label{eqnu=1.17}\\%
\mathcal{K}  &  =\operatorname*{span}\left\{  e_{-2d+1},e_{-2d+2},\dots
,e_{0}\right\}  . \label{eqnu=1.18}%
\end{align}
It follows from (\ref{eqPolynomial.23}) that the matrix for $V_{0}^{\ast}$
relative to the basis $\left\{  e_{0},e_{-1},\dots,e_{-2d+1}\right\}  $ has
the form (passing under the name ``slant-Toeplitz matrix'')%
\begin{equation}%
{\setlength{\arraycolsep}{0pt}\left( \begin{array}{cccccccccccccccccc}%
\cline{1-2}
a_0 \;& \vline&\; 0 \;&\; 0 \;& &\; 0 \;&\; 0 \;& &\; \dots\;&\; \dots
\;& &\; 0 \;&\; 0 \;& &\; 0 \;&\; 0 \;& &\; 0  \\ \cline{2-17}
a_2 \;& \vline&\; a_1 \;&\; a_0 \;& \vline&\; 0 \;&\; 0 \;& &\; \dots
\;&\; \dots\;& &\; 0 \;&\; 0 \;& &\; 0 \;&\; 0 \;& \vline&\; 0  \\ \cline{5-8}
a_4 \;& \vline&\; a_3 \;&\; a_2 \;& &\; a_1 \;&\; a_0 \;& \vline
&\; \dots\;&\; \dots\;& &\; 0 \;&\; 0 \;& &\; 0 \;&\; 0 \;& \vline&\; 0  \\
\vdots\;& \vline&\; \vdots\;&\; \vdots\;& &\; \vdots\;&\; \vdots
\;& &\; \;&\; \;& &\; \vdots\;&\; \vdots\;& &\; \vdots\;&\; \vdots
\;& \vline&\; \vdots\\ \cline{11-14}
a_{2d-4} \;& \vline&\; a_{2d-5} \;&\; a_{2d-6} \;& &\; a_{2d-7} \;&\; a_{2d-8}%
\;& &\; \dots\;&\; \dots\;& &\; a_1 \;&\; a_0 \;& \vline
&\; 0 \;&\; 0 \;& \vline&\; 0  \\ \cline{14-17}
a_{2d-2} \;& \vline&\; a_{2d-3} \;&\; a_{2d-4} \;& &\; a_{2d-5} \;&\; a_{2d-6}%
\;& &\; \dots\;&\; \dots\;& &\; a_3 \;&\; a_2 \;& &\; a_1 \;&\; a_0 \;& \vline
&\; 0  \\ \cline{1-2}\cline{17-18}
0 \;& \vline&\; a_{2d-1} \;&\; a_{2d-2} \;& &\; a_{2d-3} \;&\; a_{2d-4}%
\;& &\; \dots\;&\; \dots\;& &\; a_5 \;&\; a_4 \;& &\; a_3 \;&\; a_2 \;& \vline
&\; a_1  \\ \cline{2-5}
0 \;& \vline&\; 0 \;&\; 0 \;& \vline&\; a_{2d-1} \;&\; a_{2d-2} \;& &\; \dots
\;&\; \dots\;& &\; a_7 \;&\; a_6 \;& &\; a_5 \;&\; a_4 \;& \vline
&\; a_3  \\ \cline{5-8}
\vdots\;& \vline&\; \vdots\;&\; \vdots\;& &\; \vdots\;&\; \vdots
\;& &\; \;&\; \;& &\; \vdots\;&\; \vdots\;& &\; \vdots\;&\; \vdots
\;& \vline&\; \vdots\\
0 \;& \vline&\; 0 \;&\; 0 \;& &\; 0 \;&\; 0 \;& &\; \dots\;&\; \dots
\;& \vline&\; a_{2d-1} \;&\; a_{2d-2} \;& &\; a_{2d-3} \;&\; a_{2d-4}%
\;& \vline&\; a_{2d-5}\\ \cline{11-14}
0 \;& \vline&\; 0 \;&\; 0 \;& &\; 0 \;&\; 0 \;& &\; \dots\;&\; \dots
\;& &\; 0 \;&\; 0 \;& \vline&\; a_{2d-1} \;&\; a_{2d-2} \;& \vline
&\; a_{2d-3}\\ \cline{2-17}
0 \;& &\; 0 \;&\; 0 \;& &\; 0 \;&\; 0 \;& &\; \dots\;&\; \dots
\;& &\; 0 \;&\; 0 \;& &\; 0 \;&\; 0 \;& \vline&\; a_{2d-1}  \\ \cline{17-18}
\end{array}\right) }%
\label{eqnu=1.19}%
\end{equation}
and the matrix for $V_{1}^{\ast}$ is, by (\ref{eqnu=1.11}), obtained by using
the substitution
$a_{k}\rightarrow\left(  -1\right)  ^{k}a_{2d-1-k}$
in the matrix (\ref{eqnu=1.19}). Note that the subspace
\begin{equation}
\mathcal{K}_{0}=\operatorname*{span}\left\{  e_{-2d+2},e_{-2d+3},\dots
,e_{-1}\right\}  \label{eqnu=1.21}%
\end{equation}
is also invariant under $V_{0}^{\ast}$ and $V_{1}^{\ast}$, and thus under
$S_{0}^{\ast}$ and $S_{1}^{\ast}$, but we will see in Section
\ref{Casetheta=halfpi} below that this subspace is not always cyclic.

Let us remark that the scaling relations for
the father function $\varphi$ corresponding to (\ref{eqnu=1.10})
and the
mother function $\psi$ from (\ref{eqnu=1.11}) (both in
$L^{2}\left(\mathbb{R}\right)$) are as follows:
\begin{align}
\frac{1}{\sqrt{2}}\varphi\left(\frac{x}{2}\right)&=\sum_{k}
a_{k}\varphi\left(x-k\right) , \label{eqScalingRelation1} \\
\frac{1}{\sqrt{2}}\psi\left(\frac{x}{2}\right)&=\sum_{k}
\left(-1\right)^{k}a_{2d-1-k}\varphi\left(x-k\right) .
\label{eqScalingRelation2}
\end{align}
See also Remark \ref{RemGeometricallyClear}.

Following the
terminology in
\cite{Wel93}, we
say that $d$
is the \emph{genus},
and we now
turn to a closer
study of $d\leq 2$.

\begin{figure}
\setlength{\unitlength}{24pt}
\mbox{\begin{picture}(3,4)(-1.5,-0.7)
\put(-1.5,2.5){\makebox(3,0.5){Values of $\theta$}}
\put(-1.5,2){\makebox(3,0.5){shown:}}
\put(-0.25,1){\makebox(0.5,0.5){$\frac{\pi}{2}$}}
\put(1,-0.25){\makebox(0.5,0.5){$0$}}
\put(-0.25,-1.5){\makebox(0.5,0.5){$\frac{3\pi}{2}$}}
\put(-1.5,-0.25){\makebox(0.5,0.5){$\pi$}}
\put(-1.2812,-0.6523){\makebox(0.4,0.4){c}}
\put(-1.1998,-0.8636){\makebox(0.4,0.4){a}}
\put(0.7331,-0.8636){\makebox(0.4,0.4){b}}
\put(-1.0485,-1.0485){\makebox(0.4,0.4){d}}
\put(-1,-1){\makebox(2,2){\psfig
{figure=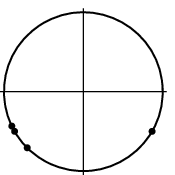,bbllx=0bp,bblly=0bp,bburx=48bp,bbury=48bp,width=48pt}}}
\end{picture}%
\hspace*{4pt}\psfig
{figure=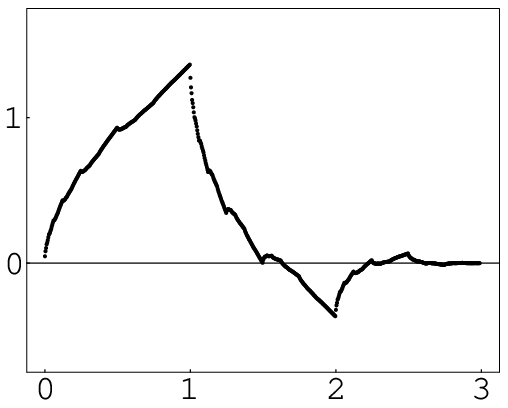,bbllx=0bp,bblly=18bp,bburx=144bp,bbury=133bp,width=120pt}%
\hspace*{24pt}\psfig
{figure=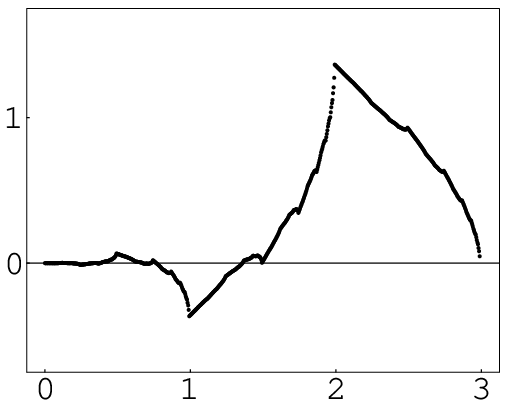,bbllx=0bp,bblly=18bp,bburx=144bp,bbury=133bp,width=120pt}%
\hspace*{20pt}}
\mbox{\hspace*{76pt}\makebox[120pt]
{\hspace*{7pt}a: $\varphi$ for $\theta=\frac{7\pi}{6}$}\hspace*{24pt}\makebox[120pt]
{\hspace*{7pt}b: $\varphi$ for $\theta=\frac{11\pi}{6}$ (or $\theta=-\frac{\pi}{6}$)}%
\hspace*{20pt}}
\mbox{\hspace*{67.5pt}\psfig
{figure=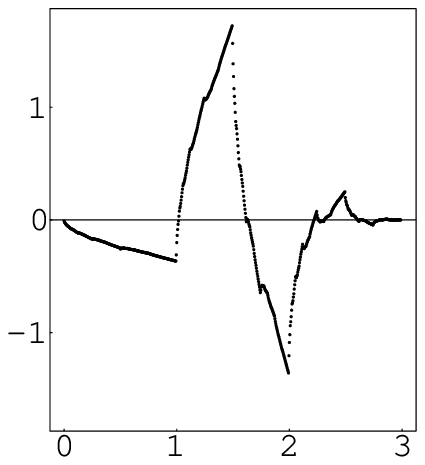,bbllx=6bp,bblly=4bp,bburx=126bp,bbury=135bp,width=129.25pt}%
\hspace*{14.75pt}\psfig
{figure=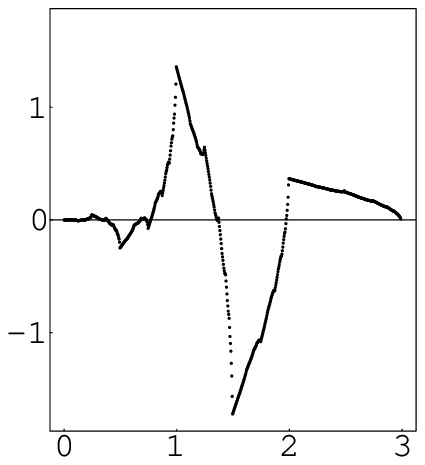,bbllx=6bp,bblly=4bp,bburx=126bp,bbury=135bp,width=129.25pt}%
\hspace*{19.25pt}}
\mbox{\hspace*{76pt}\makebox[120pt]
{\hspace*{7pt}a$^\prime$: $\psi$ for $\theta=\frac{7\pi}{6}$}\hspace*{24pt}\makebox[120pt]
{\hspace*{7pt}b$^\prime$: $\psi$ for $\theta=\frac{11\pi}{6}$ (or $\theta=-\frac{\pi}{6}$)}%
\hspace*{20pt}}
\mbox{\hspace*{76pt}\psfig
{figure=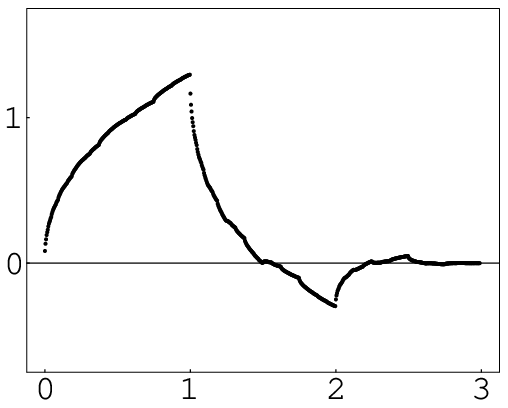,bbllx=0bp,bblly=18bp,bburx=144bp,bbury=136bp,width=120pt}%
\hspace*{24pt}\psfig
{figure=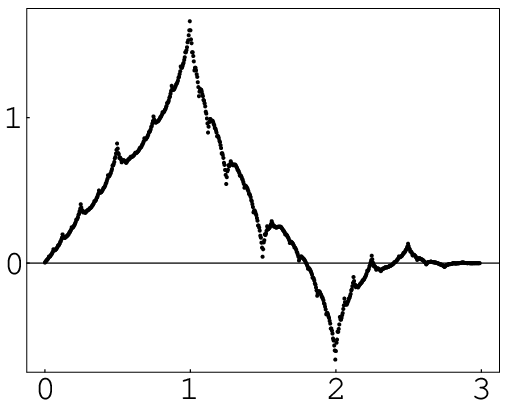,bbllx=0bp,bblly=18bp,bburx=144bp,bbury=136bp,width=120pt}%
\hspace*{20pt}}
\mbox{\hspace*{76pt}\makebox[120pt]
{\hspace*{7pt}c: $\varphi$ for $\theta=\frac{8\pi}{7}$}\hspace*{24pt}\makebox[120pt]
{\hspace*{7pt}d: $\varphi$ for $\theta=\frac{5\pi}{4}$}\hspace*{20pt}}
\mbox{\hspace*{67.5pt}\psfig
{figure=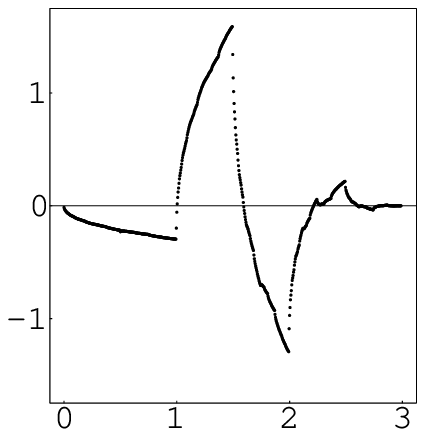,bbllx=2bp,bblly=4bp,bburx=122bp,bbury=127bp,width=129.25pt}%
\hspace*{14.75pt}\psfig
{figure=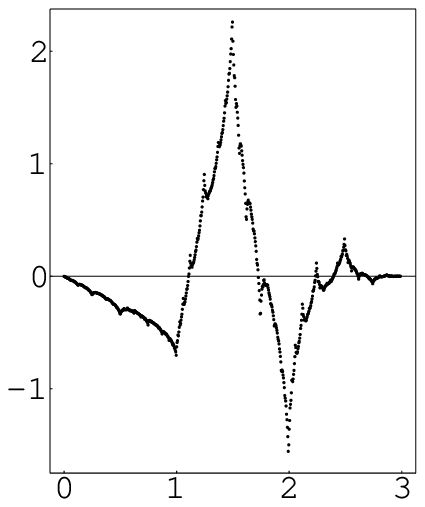,bbllx=12bp,bblly=4bp,bburx=132bp,bbury=147bp,width=129.25pt}%
\hspace*{19.25pt}}
\mbox{\hspace*{76pt}\makebox[120pt]
{\hspace*{7pt}c$^\prime$: $\psi$ for $\theta=\frac{8\pi}{7}$}\hspace*{24pt}\makebox[120pt]
{\hspace*{7pt}d$^\prime$: $\psi$ for $\theta=\frac{5\pi}{4}$}\hspace*{20pt}}
\caption{Father ($\varphi$) and mother ($\psi$) functions for $\theta$ near $\frac{7\pi}%
{6}$ and $-\frac{\pi}{6}$: Continuous cases (Case ``a'' = Daubechies wavelet)}
\label{FatherMotherNear}
\end{figure}

\subsubsection{\label{Casenu=1N=2d=1}The case with dimension $\nu=1$, scale
$N=2$, and genus $d=1$}

In this case,
the second condition of (\ref{eqnu=1.12and13}) is vacuous, and the only solution of
(\ref{eqnu=1.12and13}) and (\ref{eqnu=1.14}) is $a_{0}=a_{1}=\frac{1}{\sqrt{2}}$,
so%
\begin{equation}
m_{0}\left(  z\right)  =\left( 1+z\right) /\sqrt{2},\qquad m_{1}\left(  z\right)
=\left( 1-z\right) /\sqrt{2}, \label{eqnu=1.22}%
\end{equation}
which is exactly the Haar wavelet
(Figure \ref{HaarMother}). The
representation splits into the direct sum of the two inequivalent irreducible
representations in (\ref{eqIntro.23}) and (\ref{eqIntro.24}), and the
restriction of each of these representations to $\operatorname*{UHF}%
\nolimits_{2}$ is still irreducible by \cite[Proposition 8.1]{BrJo96b}. This
can also be checked directly: in this case,%
\begin{equation}
V_{0}^{\ast}=%
\begin{pmatrix}
\frac{1}{\sqrt{2}} & 0\\
0 & \frac{1}{\sqrt{2}}%
\end{pmatrix}
,\qquad V_{1}^{\ast}=%
\begin{pmatrix}
\frac{1}{\sqrt{2}} & 0\\
0 & -\frac{1}{\sqrt{2}}%
\end{pmatrix}
. \label{eqnu=1.23}%
\end{equation}
Thus%
\begin{equation}
\bs%
\begin{pmatrix}
a & b\\
c & d
\end{pmatrix}
=\sum_{i=0}^{1}V_{i}^{{}}%
\begin{pmatrix}
a & b\\
c & d
\end{pmatrix}
V_{i}^{\ast}=%
\begin{pmatrix}
a & 0\\
0 & d
\end{pmatrix}
, \label{eqnu=1.24}%
\end{equation}
so $\mathcal{B}\left(  \mathcal{K}\right)  ^{\bs}$ is the $\ast$-algebra of
all diagonal $2\times2$ matrices. Thus the representation splits into the
direct sum of two representations with the one-dimensional $S_{i}^{\ast}%
$-invariant subspaces $\mathbb{C}%
\mkern2mu%
e_{0}$ and $\mathbb{C}%
\mkern2mu%
e_{-1}$. The corresponding maps $\bs$ on the one-dimensional subspaces are
both equal to the identity, thus they are ergodic with peripheral spectrum
$1$, and $\operatorname*{UHF}\nolimits_{2}$ is dense by Corollary
\ref{CorPolynomial.3}. Note that the states on $\mathcal{O}_{2}$ corresponding
to $e_{0}$ and $e_{-1}$ are the Cuntz states
(see \cite{Cun77,
Eva80,
BJKW97})
\begin{equation}
\omega_{0}\left(  S_{I}^{{}}S_{J}^{\ast}\right)   =2^{-\frac{\left|
I\right|  +\left|  J\right|  }{2}},\qquad
\omega_{-1}\left(  S_{I}^{{}}S_{J}^{\ast}\right)   =\left(  -1\right)
^{\left|  I\right|  +\left|  J\right|  }2^{-\frac{\left|  I\right|  +\left|
J\right|  }{2}}. \label{eqnu=1.25and26}%
\end{equation}

\subsubsection{\label{Casenu=1N=2d=2}The case with dimension $\nu=1$, scale
$N=2$, and genus $d=2$}

We now display the one-parameter family (with two singular points) of mutually
inequivalent irreducible representations of $\mathcal{O}_{2}$ mentioned in the
Introduction, and we relate the representation-theoretic behavior to the
corresponding properties of the associated family of wavelets on $L^{2}\left(
\mathbb{R}\right)  $. In this case the algebraic variety defined by
(\ref{eqnu=1.12and13})--(\ref{eqnu=1.14}) is actually the circle, and may be
defined by the following parametrization:%
\begin{equation}
\begin{aligned} a_{0} & =\frac{1}{2\sqrt{2}}\left( 1-\cos\theta+\sin
\theta\right) ,&\qquad a_{1} & =\frac{1}{2\sqrt{2}}\left( 1-\cos\theta-\sin
\theta\right) ,\\ a_{2} & =\frac{1}{2\sqrt{2}}\left( 1+\cos\theta-\sin
\theta\right) ,&\qquad a_{3} & =\frac{1}{2\sqrt{2}}\left( 1+\cos\theta+\sin
\theta\right) ; \end{aligned} \label{eqnu=1.27}%
\end{equation}
see
\cite{Pol89,Pol90}, and also
\cite{Wel93,HRW92}.
Let us give a simple argument for
this parametrization: View $a=\left(a_0,a_1,a_2,a_3\right)$
as a function on
the cyclic group of order $4$,
$\mathbb{Z}_{4}$, and consider the
Fourier transform on $\mathbb{Z}_{4}$:
\[
\hat{a}\left(n\right)=\frac{1}{2}\sum_{m=0}^{3}i^{nm}a\left(m\right),  \qquad
a\left(m\right)=\frac{1}{2}\sum_{n=0}^{3}i^{-nm}\hat{a}\left(n\right).
\]
We have the usual formulae
\[
\sum_{m}\overline{a\left(m\right)}b\left(m\right)
=\sum_{n}\overline{\hat{a}\left(n\right)}\hat{b}\left(n\right),  \qquad
\widehat{a\left(\,\cdot\,+k\right)}\left(n\right)
=i^{-nk}\hat{a}\left(n\right),  
\]
and thus
\[
\sum_{m}\overline{a\left(m+2\right)}a\left(m\right)
=\sum_{n}\left(-1\right)^{n}\overline{\hat{a}\left(n\right)}
\hat{a}\left(n\right).
\]
Also
$\Hat{\Bar{a}}\left(n\right)=\overline{\hat{a}\left(-n\right)}$.
The relations (\ref{eqnu=1.4and5}) and (\ref{eqnu=1.7}), 
together with reality of $a$, take the form
\begin{gather*}
\sum_{n}a_{n}=\sqrt{2},  \qquad
\sum_{n}a_{n}^{2}=\sum_{n}\bar{a}_{n}a_{n}=1,  \qquad
a_{n}=\bar{a}_{n},  \\
a_{0}a_{2}+a_{1}a_{3}=0
\Longleftrightarrow\sum_{m}\overline{a\left(m+2\right)}a\left(m\right)=0,
\end{gather*}
and hence
\begin{gather*}
\hat{a}\left(0\right)=\frac{1}{\sqrt{2}},  \qquad
\sum_{n}\left|\hat{a}\left(n\right)\right|^{2}=1,  \\
\overline{\hat{a}\left(-n\right)}=a\left(n\right)
\Longleftrightarrow\hat{a}\left(0\right),\hat{a}\left(2\right)
\text{ are real and }\hat{a}\left(3\right)=\overline{\hat{a}\left(1\right)}, \\
\sum_{n}\left(-1\right)^{n}\overline{\bar{a}\left(n\right)}
\hat{a}\left(n\right)=0.
\end{gather*}
Introducing $c=a\left(2\right)=\bar{c}$ and $b=\hat{a}\left(1\right)$
we thus have
$c^{2}+2\left|b\right|^{2}=1/2$,
$c^{2}-2\left|b\right|^{2}=-1/2$,
and hence
$c=0$, $\left|b\right|=1/2$.
Putting $b=\frac{1}{2}e^{i\varphi}$, the relations
for $a$ are thus equivalent to
\[
\left(\hat{a}\left(0\right),\hat{a}\left(1\right),\hat{a}\left(2\right),\hat{a}\left(3\right)\right)
=\left(\frac{1}{\sqrt{2}},\frac{e^{i\varphi}}{2},0,\frac{e^{-i\varphi}}{2}\right).
\]
Applying the inverse Fourier transform
to this, we obtain
\begin{align*}
a_{0} & =\frac{1}{2\sqrt{2}}\left( 1+\sqrt{2}\cos\varphi\right) ,&
a_{1} & =\frac{1}{2\sqrt{2}}\left( 1+\sqrt{2}\sin\varphi\right) ,\\
a_{2} & =\frac{1}{2\sqrt{2}}\left( 1-\sqrt{2}\cos\varphi\right) ,&
a_{3} & =\frac{1}{2\sqrt{2}}\left( 1-\sqrt{2}\sin\varphi\right) .
\end{align*}
Substituting $\varphi=\theta+\frac{5\pi}{4}$ here, we
obtain (\ref{eqnu=1.27}).

Returning to the representation, the
operators $V_{k}^{*}$ from (\ref{eqnu=1.19}) in this
case have the form:
\begin{equation}
V_{0}^{\ast}  =%
{\setlength{\arraycolsep}{0pt}\left( \begin{array}{cccccc}
a_0 \;& &\; 0 \;&\; 0 \;& &\; 0  \\ \cline{2-5}
a_2 \;& \vline&\; a_1 \;&\; a_0 \;& \vline&\; 0  \\
0 \;& \vline&\; a_3 \;&\; a_2 \;& \vline&\; a_1  \\ \cline{2-5}
0 \;& &\; 0 \;&\; 0 \;& &\; a_3
\end{array}\right) }%
\text{\quad and\quad }
V_{1}^{\ast}    =%
{\setlength{\arraycolsep}{0pt}\left( \begin{array}{cccccc}
a_3 \;& &\; 0 \;&\; 0 \;& &\; 0  \\ \cline{2-5}
a_1 \;& \vline&\; -a_2 \;&\; a_3 \;& \vline&\; 0  \\
0 \;& \vline&\; -a_0 \;&\; a_1 \;& \vline&\; -a_2  \\ \cline{2-5}
0 \;& &\; 0 \;&\; 0 \;& &\; -a_0
\end{array}\right) }%
. \label{eqnu=1.28and29}%
\end{equation}
If one replaces the angle variable $\theta$ with $\varphi$, and calls the
corresponding coefficients $b_{0},\dots,b_{3}$, and the corresponding matrices
$W_{0}^{\ast}$, $W_{1}^{\ast}$, the corresponding map $\br\colon
M_{4}\rightarrow M_{4}$ given by (\ref{eqFinitely.17}),%
\begin{equation}
\br\left(  A\right)  =\sum_{i=0}^{1}W_{i}^{{}}AV_{i}^{\ast}, \label{eqnu=1.30}%
\end{equation}
is defined by a $16\times16$ matrix relative to the basis
\begin{equation}
e_{0,0},\,e_{0,-1},\,\dots,\,e_{0,-3},\,e_{-1,0},\,e_{-1,-1},\,\dots
,\,e_{-3,-3} \label{eqBasis16}%
\end{equation}
of $M_{4}$. This $16\times16$ matrix has the form%
\begin{equation}%
{\setlength{\arraycolsep}{0pt}\left( \begin{array}{cccccc}
A_0 \;& &\; A_2 \;&\; 0 \;& &\; 0  \\ \cline{2-5}
0 \;& \vline&\; A_1 \;&\; A_3 \;& \vline&\; 0  \\
0 \;& \vline&\; A_0 \;&\; A_2 \;& \vline&\; 0  \\ \cline{2-5}
0 \;& &\; 0 \;&\; A_1 \;& &\; A_3
\end{array}\right) }%
, \label{eqnu=1.31}%
\end{equation}
where the $4\times4$ matrices $A_{i}$ are given by%
\begin{equation}
\begin{aligned} A_{0} & =b_{0}V_{0}+b_{3}V_{1},&\qquad
A_{1} & =b_{1}V_{0}%
-b_{2}V_{1},\\ A_{2} & =b_{2}V_{0}+b_{1}V_{1},&\qquad
A_{3} & =b_{3}V_{0}%
-b_{0}V_{1}. \end{aligned} \label{eqnu=1.32}%
\end{equation}
Thus one can compute the eigenvalues of $\br$ by computing the eigenvalues of
the matrices $A_{0}$, $A_{3}$, and $\left(
\begin{smallmatrix}
A_{1} & A_{3}\\
A_{0} & A_{2}%
\end{smallmatrix}
\right)  $. If $\varphi=\theta$ the result is (we call $\br=\bs$ in this case
as usual)%
\begin{equation}%
\newcommand{\eigenstrut}{\vphantom{$\displaystyle
\left( \frac{1+\sin\theta}{2}\right) $}}
\begin{tabular}{r|cccccc}
\eigenstrut Eigenvalue of $\bs$ &
$1$ &
$0$ &
$\displaystyle\frac{\cos\theta}{2}$ &
$\displaystyle-\frac{\cos\theta}{2}$ &
$\displaystyle\frac{1+\sin\theta}{2}$ &
$-\sin\theta$ \\ \hline
\eigenstrut Multiplicity  
& $1$  
& $8$  
& $2$  
& $2$  
& $2$  
& $1$.
\end{tabular}%
\label{eqnu=1.33}%
\end{equation}
Hence, the dimension of the eigenspace $\left\{  A\mid\bs\left(  A\right)
=A\right\}  $ is%
\begin{equation}%
\begin{cases}
3 & \text{if }\theta=\frac\pi2 \\
2 & \text{if }\theta=\frac{3\pi}2 \\
1 & \text{otherwise.}
\end{cases}%
\label{eqnu=1.34}%
\end{equation}
These numbers are then the dimensions of the commutants of the corresponding
representations. Since the only $C^{\ast}$-algebras of dimensions $1$, $2$,
$3$ are $\mathbb{C}$, $\mathbb{C}%
\mkern1mu%
^{2}$, $\mathbb{C}%
\mkern1mu%
^{3}$, it follows that the representation of $\mathcal{O}_{2}$ splits into $2$
inequivalent irreducible representations if $\theta=\frac{3\pi}{2}$, into $3$
inequivalent irreducible representations when $\theta=\frac{\pi}{2}$, and the
representation is irreducible for all other $\theta$. We note that the
peripheral spectrum of $\bs$ is nontrivial only if $\theta=\frac{\pi}{2}$,
when $-1$ is an eigenvalue of multiplicity $1$. Thus the representations for
generic $\theta\notin\left\{  \frac{\pi}{2},\frac{3\pi}{2}\right\}  $ also
have irreducible restriction to $\operatorname*{UHF}\nolimits_{2}$. Finally,
if one considers the case $\theta\neq\varphi$, one can compute that $1$ is an
eigenvalue for $\br$ if and only if $\left\{  \theta,\varphi\right\}
=\left\{  \frac{\pi}{2},\frac{3\pi}{2}\right\}  $, and the dimension of the
corresponding eigenspace is then $2$. We recall from (\ref{eqFinitely.17})
that solutions $A\neq0$ to $\br\left(  A\right)  =A$ correspond by lifting to
operators on $L^{2}\left(  \mathbb{T}\right)  $ which intertwine the two
associated $\mathcal{O}_{2}$-representations $\bp^{\left(  \theta\right)  }$
and $\bp^{\left(  \varphi\right)  }$ for $\theta$ and $\varphi$, respectively.
Hence the representations for generic points $\theta\notin\left\{  \frac{\pi
}{2},\frac{3\pi}{2}\right\}  $ on the circle are all mutually disjoint by
(\ref{eqFinitely.16})--(\ref{eqFinitely.18}), but if $\left\{  \theta
,\varphi\right\}  =\left\{  \frac{\pi}{2},\frac{3\pi}{2}\right\}  $, the
intertwiner space is $2$-dimensional. See Section \ref{Intertwining} for more
details on the latter.

\begin{figure}
\mbox{\psfig
{figure=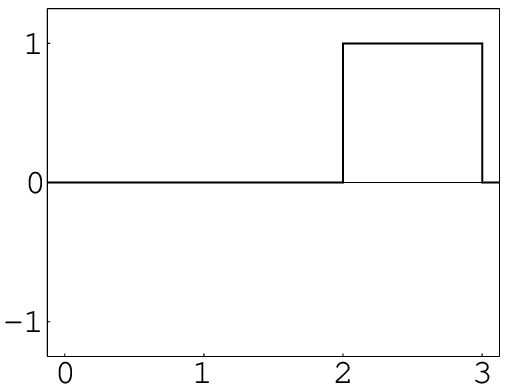,bbllx=0bp,bblly=20bp,bburx=144bp,bbury=145bp,width=120pt}%
\psfig
{figure=FatherHaar.eps,bbllx=0bp,bblly=20bp,bburx=144bp,bbury=145bp,width=120pt}%
\psfig
{figure=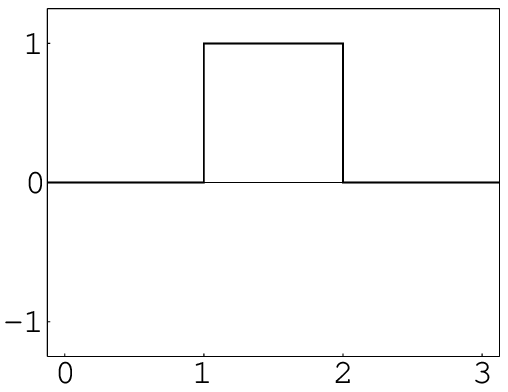,bbllx=0bp,bblly=20bp,bburx=144bp,bbury=145bp,width=120pt}}
\mbox{\makebox[120pt]
{\hspace*{10pt}$\varphi$ for $\theta=0$}\makebox[120pt]
{\hspace*{10pt}$\varphi$ for $\theta=\pi$}\begin{minipage}
[t]{120pt}\makebox[120pt]
{\hspace*{10pt}$\varphi$ for $\theta=\frac{3\pi}{2}$}
\makebox[120pt]{\hspace*{10pt}(cf.\ (\ref{eqnu=1.36}))}\end{minipage}}
\mbox{\psfig
{figure=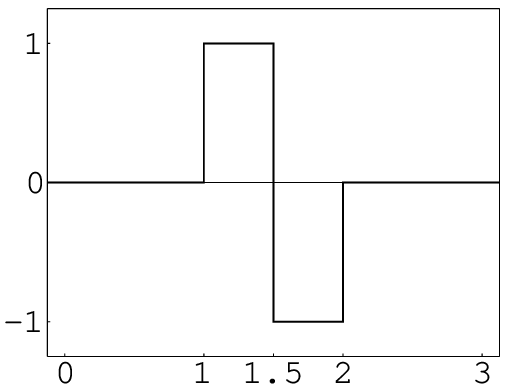,bbllx=0bp,bblly=20bp,bburx=144bp,bbury=145bp,width=120pt}%
\psfig
{figure=Mother000.eps,bbllx=0bp,bblly=20bp,bburx=144bp,bbury=145bp,width=120pt}%
\psfig
{figure=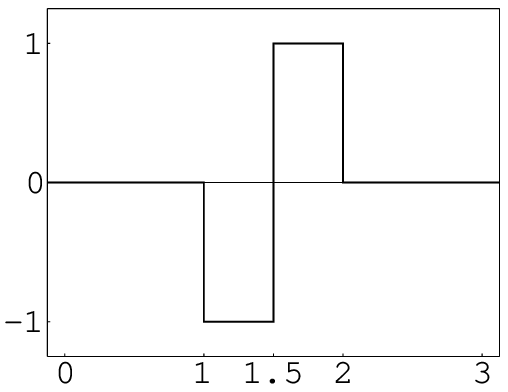,bbllx=0bp,bblly=20bp,bburx=144bp,bbury=145bp,width=120pt}}
\mbox{\begin{minipage}[t]{120pt}\makebox[120pt]
{\hspace*{10pt}$\psi$ for $\theta=0$}
\makebox[120pt]
{\hspace*{10pt}(cf.\ (\ref{eqnu=1.68and69}))}\end{minipage}\begin{minipage}
[t]{120pt}\makebox[120pt]
{\hspace*{10pt}$\psi$ for $\theta=\pi$}
\makebox[120pt]
{\hspace*{10pt}(cf.\ (\ref{eqnu=1.68and69}))}\end{minipage}\begin{minipage}
[t]{120pt}\makebox[120pt]
{\hspace*{10pt}$\psi$ for $\theta=\frac{3\pi}{2}$}
\makebox[120pt]
{\hspace*{10pt}(cf.\ (\ref{eqnu=1.37}))}\end{minipage}}
\caption{Father ($\varphi$) and mother ($\psi$) functions
for $\theta$ equal to multiples of $\frac
{\pi
}{2}$: The symmetry $\psi(3-x)=-\psi(x)$}
\label{FatherMothertheta=khalfpi}
\end{figure}

\begin{figure}
\mbox{\psfig
{figure=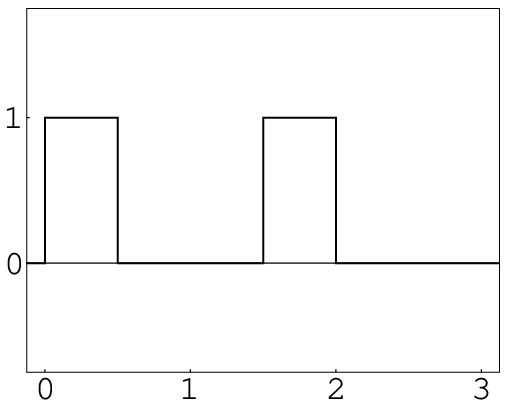,bbllx=0bp,bblly=18bp,bburx=144bp,bbury=150bp,width=120pt}%
\psfig
{figure=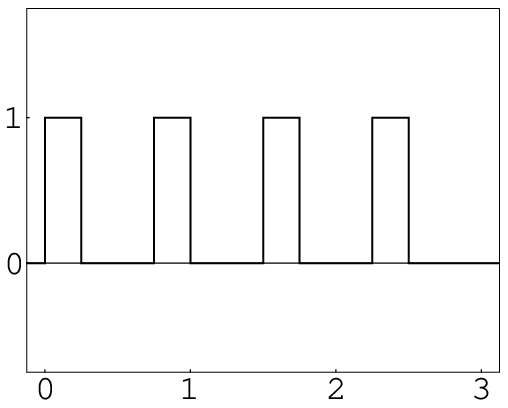,bbllx=0bp,bblly=18bp,bburx=144bp,bbury=150bp,width=120pt}%
\psfig
{figure=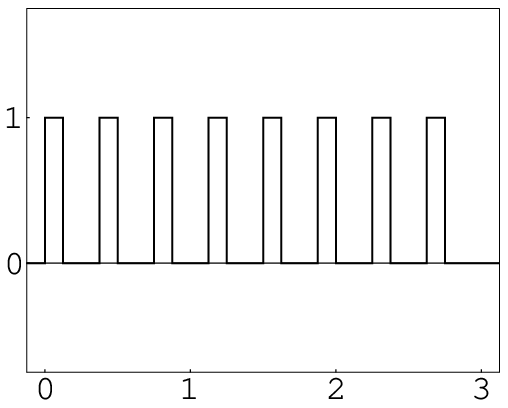,bbllx=0bp,bblly=18bp,bburx=144bp,bbury=150bp,width=120pt}}
\mbox{\makebox[120pt]
{\hspace*{8pt}First approximation to $\varphi$}\makebox[120pt]
{\hspace*{8pt}Second approximation}\makebox[120pt]
{\hspace*{8pt}Third approximation}}
\mbox{\psfig
{figure=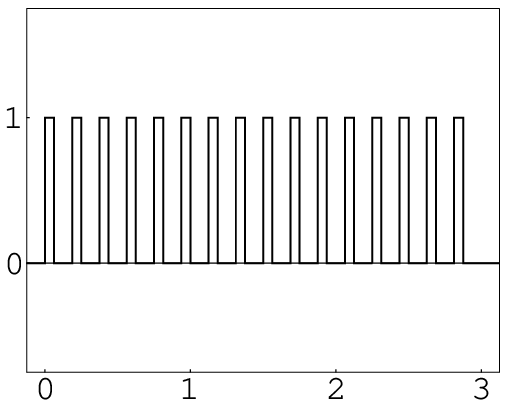,bbllx=0bp,bblly=18bp,bburx=144bp,bbury=150bp,width=120pt}%
\psfig
{figure=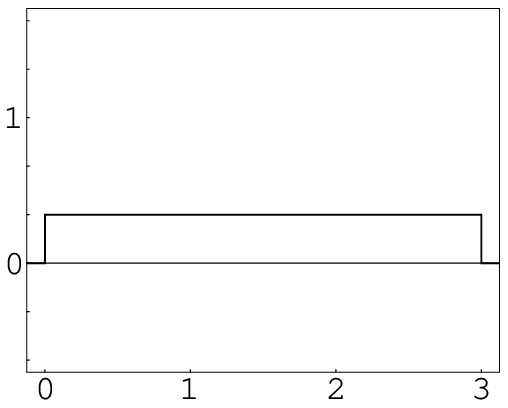,bbllx=0bp,bblly=18bp,bburx=144bp,bbury=150bp,width=120pt}%
\psfig
{figure=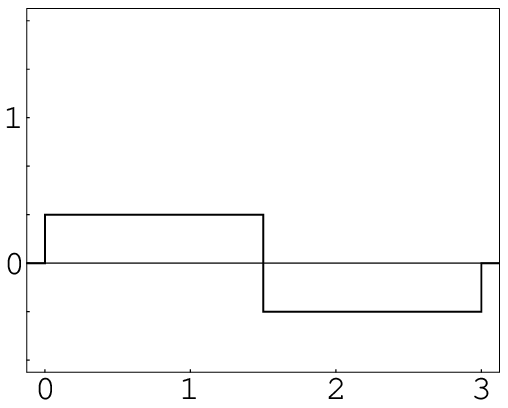,bbllx=0bp,bblly=18bp,bburx=144bp,bbury=150bp,width=120pt}}
\mbox{\makebox[120pt]
{\hspace*{8pt}Fourth approximation}\begin{minipage}
[t]{120pt}\makebox[120pt]
{\hspace*{8pt}Weak $L^{2}$ limit: $\varphi$ for $\theta=\frac{\pi}{2}$}
\makebox[120pt]
{\hspace*{8pt}(cf.\ (\ref{eqnu=1.47}))}\end{minipage}\begin{minipage}
[t]{120pt}\makebox[120pt]
{\hspace*{8pt}$\psi$ for $\theta=\frac{\pi}{2}$}
\makebox[120pt]
{\hspace*{8pt}(cf.\ (\ref{eqnu=1.48}))}\end{minipage}}
\renewcommand{\frogcapReptileGenealogy}{: See discussion in
Section \ref{Continuity}}
\caption{Father ($\varphi$) and mother ($\psi$) functions
for $\theta=\frac{\pi
}{2}$, with cascade-algorithm approximations of father function
$\varphi$\protect\frogcapReptileGenealogy}
\label{ReptileGenealogy}
\end{figure}

A second immediate observation on (\ref{eqnu=1.27}) is that at the four points
$\theta=0,\frac{\pi}{2},\pi,\frac{3\pi}{2}$, we have two of the four
coefficients vanishing with different pairs in the four different cases, so
those four cases are closely connected to four modified Haar wavelets,
illustrated in Figures \ref{ReptileGenealogy} and
\ref{FatherMothertheta=khalfpi}.
A more subtle fact, to be described below, is that it is only
the two cases $\theta=\frac{\pi}{2}$ and $\theta=\frac{3\pi}{2}$ on the
symmetry axis where the corresponding $\mathcal{O}_{2}$-representation on
$L^{2}\left(  \mathbb{T}\right)  $ fails to be irreducible. The case
$\theta=\frac{\pi}{2}$ is degenerate in a sense illustrated in Figure
\ref{ReptileGenealogy}. We will relate the resulting degenerate decomposition
at $\theta=\frac{\pi}{2}$ of the subalgebra $\operatorname*{UHF}%
\nolimits_{2}\subset\mathcal{O}_{2}$ to the wavelet properties.

Let us now consider the two exceptional points $\theta=\frac{\pi}{2}$ and
$\theta=\frac{3\pi}{2}$ separately.

\subsubsubsection{\label{Casetheta=3halfpi}The case $\theta=\frac{3\pi}{2}$}

When $\theta=\frac{3\pi}{2}$,
$a_{0}=a_{3}=0$, $a_{1}=a_{2}=1/\sqrt{2}$,
so
\begin{align}
m_{0}\left(  z\right)   &  =\left(  z+z^{2}\right) /\sqrt{2} , &
\varphi\left(  x/2\right)   & =\varphi\left(  x-1\right)
+\varphi\left(  x-2\right)  ,\label{eqnu=1.36}\\
m_{1}\left(  z\right)   &  =\left(  -z+z^{2}\right) /\sqrt{2} , &
\psi\left(  x/2\right)   & =-\varphi\left(  x-1\right)  +\varphi
\left(  x-2\right)  , \label{eqnu=1.37}%
\end{align}
with the scaling relations indicated for the father function $\varphi$, and
the mother function $\psi$, respectively; see Figure
\ref{FatherMothertheta=khalfpi}.%

This is a simple transform of the Haar wavelet
(Figure \ref{HaarMother}), and the representation theory becomes similar: defining
$S_{i}$ by (\ref{eqIntro.12}) and transforming the representation by $\frac
{1}{\sqrt{2}}\left(
\begin{smallmatrix}
1 & 1\\
-1 & 1
\end{smallmatrix}
\right)  \in\mathrm{U}\left(  2\right)  $, i.e.,%
\begin{equation}
T_{0} =\left( S_{0}-S_{1}\right
) /\sqrt{2},\qquad T_{1} =\left( S_{0}+S_{1}\right) /\sqrt{2}, 
\label{eqnu=1.40}%
\end{equation}
we obtain%
\begin{equation}
T_{0}\xi\left( z\right) =z\xi\left( z^{2}\right
) ,\qquad T_{1}\xi\left( z\right) =z^{2}\xi\left( z^{2}\right) .
\label{eqnu=1.41}%
\end{equation}
By the computation in \cite[eqs.~(8.1)--(8.2)]{BrJo96b}, if $U$ is the
unitary operator given by multiplication by $z^{-1}$, then%
\begin{equation}
U^{\ast}T_{0}U\xi\left( z\right) =\xi\left( z^{2}
\right) ,\qquad U^{\ast}T_{1}U\xi\left( z\right) =z\xi\left( z^{2}
\right) .
\label{eqnu=1.42}%
\end{equation}
By \cite[Proposition 8.1]{BrJo96b}, $L^{2}\left(  \mathbb{T}\right)  $ splits
into the two irreducible subspaces spanned by%
$\left\{  1,z,z^{2},\dots\right\}$
and
$\left\{  z^{-1},z^{-2},\dots\right\}  $.
Applying $U$ to these, we obtain the two irreducible invariant subspaces
corresponding to the original representation%
\begin{equation}
\overline{\operatorname*{span}}\left\{  z^{-1},1,z,z^{2},\dots\right\}
\text{\quad and\quad }
\overline{\operatorname*{span}}\left\{  z^{-2},z^{-3},\dots\right\} 
\label{eqnu=1.43and44}%
\end{equation}
(overbar for closure).
We see that the projection $P$ onto the overlapping four-dimensional
$S_{i}^{\ast}$-invariant subspace $\mathcal{K}=\operatorname*{span}\left\{
1,z^{-1},z^{-2},z^{-3}\right\}  $ commutes with the projection onto the first
two subspaces. The respective products of $P$ by these projections are%
\begin{equation}%
\begin{pmatrix}
1 & 0 & 0 & 0\\
0 & 1 & 0 & 0\\
0 & 0 & 0 & 0\\
0 & 0 & 0 & 0
\end{pmatrix}
\text{\quad and\quad}%
\begin{pmatrix}
0 & 0 & 0 & 0\\
0 & 0 & 0 & 0\\
0 & 0 & 1 & 0\\
0 & 0 & 0 & 1
\end{pmatrix}
, \label{eqnu=1.45}%
\end{equation}
and these two matrices span exactly the eigenspace of $\bs$ corresponding to
eigenvalue $1$. Also, each of the two subrepresentations has irreducible
restriction to $\operatorname*{UHF}\nolimits_{2}$, confirming the fact that
the peripheral spectrum of $\bs$ consists of $1$ alone.

\subsubsubsection{\label{Casetheta=halfpi}The case $\theta=\frac{\pi}{2}$}

When $\theta=\frac{\pi}{2}$,%
\begin{equation}
a_{0}=a_{3}=1/\sqrt{2},\qquad a_{1}=a_{2}=0, \label{eqnu=1.46}%
\end{equation}
so the associated low/high-pass filters and scaling relations are:
\begin{align}
m_{0}\left(  z\right)   &  =\left(  1+z^{3}\right) /\sqrt{2}, &
\varphi\left(  x/2\right)   & =\varphi\left(  x\right)
+\varphi\left(  x-3\right)  
,\label{eqnu=1.47}\\
m_{1}\left(  z\right)   &  =\left(  1-z^{3}\right)  /\sqrt{2}, &
\psi\left(  x/2\right)   & =\varphi\left(  x\right)  -\varphi
\left(  x-3\right) .
\label{eqnu=1.48}%
\end{align}
See Figure \ref{ReptileGenealogy} for the graphs of the corresponding
$\varphi$ and $\psi$. Applying the unitary $\frac{1}{\sqrt{2}}\left(
\begin{smallmatrix}
1 & 1\\
1 & -1
\end{smallmatrix}
\right)  \in\mathrm{U}\left(  2\right)  $ to this representation, we transform
it into the representation with $m_{0}\left(  z\right)  =1$, $m_{1}\left(
z\right)  =z^{3}$. We have already noted in (\ref{eqnu=1.34}) that the fixed
point set of $\bs$ is three-dimensional in this case, and indeed, by
\cite[Proposition 8.2]{BrJo96b}, this representation decomposes into $3$
mutually disjoint irreducible representations given by restriction to the $3$
subspaces%
\begin{equation}
\begin{gathered}
\overline{\operatorname*{span}}\left\{ z^{3n}\mid
n=0,1,2,\dots\right\} ,\qquad \overline{\operatorname*{span}}\left
\{ z^{3n}\mid n=-1,-2,\dots\right\} ,\\ \overline{\operatorname*{span}
}\left\{ z^{k}\mid k\text{ not divisible by }3\right\} .
\end{gathered}
\label{eqnu=1.49}%
\end{equation}
The restriction to $\operatorname*{UHF}\nolimits_{2}$ is still irreducible on
the first two subspaces, while it decomposes into the two irreducible
subrepresentations on%
\begin{equation}
\overline{\operatorname*{span}}\left\{ z^{3k+1}\mid
k\in\mathbb{Z}\right\} ,\qquad \overline{\operatorname*{span}}\left
\{ z^{3k+2}\mid k\in\mathbb{Z}\right\} 
\label{eqnu=1.50}%
\end{equation}
on the third subspace. Again the projection onto each of these subspaces
commutes with $P$, and hence the eigenspace of $\bs$ corresponding to
eigenvalue $1$ is spanned by the three projections%
\begin{equation}%
\begin{pmatrix}
1 & 0 & 0 & 0\\
0 & 0 & 0 & 0\\
0 & 0 & 0 & 0\\
0 & 0 & 0 & 0
\end{pmatrix}
,\quad%
\begin{pmatrix}
0 & 0 & 0 & 0\\
0 & 0 & 0 & 0\\
0 & 0 & 0 & 0\\
0 & 0 & 0 & 1
\end{pmatrix}
\text{,\quad and\quad}%
\begin{pmatrix}
0 & 0 & 0 & 0\\
0 & 1 & 0 & 0\\
0 & 0 & 1 & 0\\
0 & 0 & 0 & 0
\end{pmatrix}
, \label{eqnu=1.51}%
\end{equation}
respectively, confirming that the eigenvalue $1$ has multiplicity $3$ in this
case. Furthermore, if $U$ is the unitary operator (\ref{eqFinitely.26}) on
$\overline{\operatorname*{span}}\left\{  z^{k}\mid k\text{ not divisible by
}3\right\}  $ that implements the gauge automorphism $\tau_{-1}$ there, we
have $UT_{i}=-T_{i}U$. Hence%
\begin{equation}
U\left(  \xi\left(  z^{2}\right)  \right)   =-\left(  U\xi\right)  \left(
z^{2}\right) 
\text{\quad and\quad }
U\left(  z^{3}\xi\left(  z^{2}\right)  \right)   =-z^{3}\left(
U\xi\right)  \left(  z^{2}\right)  
\label{eqnu=1.52and53}%
\end{equation}
if $\xi$ is in this subspace. This unitary $U$ from (\ref{eqFinitely.26}) has
to fix the two subspaces $\overline{\operatorname*{span}}\left\{  z^{3k+1}\mid
k\in\mathbb{Z}\right\}  $ and $\overline{\operatorname*{span}}\left\{
z^{3k+2}\mid k\in\mathbb{Z}\right\}  $, and $U^{2}=\openone$, hence it is
clear that%
\begin{equation}
PUP=\pm%
\begin{pmatrix}
0 & 0 & 0 & 0\\
0 & 1 & 0 & 0\\
0 & 0 & -1 & 0\\
0 & 0 & 0 & 0
\end{pmatrix}
. \label{eqnu=1.54}%
\end{equation}
It is easily verified directly that this is the eigenvector of $\bs$
corresponding to eigenvalue $-1$. This means that the group from
(\ref{eqFinitely.21})--(\ref{eqFinitely.23}) in this case is $\mathbb{Z}_{2}$,
if $\mathcal{K}$ is taken to be $\operatorname*{span}\left\{  e_{-1}%
,e_{-2}\right\}  $ and $\varphi$ the trace state on $\mathcal{B}\left(
\mathcal{K}\right)  $.

In conclusion we note that this $\mathcal{O}_{2}$-representation $\bp^{\left(
\theta\right)  }$, $\theta=\frac{\pi}{2}$, as well as its restriction to
$\operatorname*{UHF}\nolimits_{2}$ has a decomposition into irreducibles which
sets it apart from the other representations when $\theta\neq\frac{\pi}{2}$.
We
will see in the beginning of Section \ref{Symmetry}
that if $\theta=0$, $\pi$, or $\frac{3\pi}{2}$, then the wavelet is
still of Haar type, i.e., $\varphi$ is of the form $\varphi=\chi_{I}$ where
$I$ is an interval of unit length. The position of the interval $I$ varies
(see Figure \ref{FatherMothertheta=khalfpi})
in
the three cases $\theta=0$, $\pi$, or $\frac{3\pi}{2}$, while the
mother function $\psi^{\left(  0\right)  }$ is common for two of them,
$\theta=0,\pi$, and
$\psi^{\left(  \frac{3\pi}{2}\right)  }=-\psi^{\left(  0\right)  }
=-\psi^{\left(  \pi\right)  }$.
All three satisfy $\psi\left(  3-x\right)  =-\psi\left(  x\right)  $. But, if
$\theta=\frac{\pi}{2}$, then the nature of $\varphi$ is somewhat different.
{}From (\ref{eqnu=1.47}), we see that $\varphi\left(  \frac{x}{2}\right)
=\varphi\left(  x\right)  +\varphi\left(  x-3\right)  $; and by \cite{GrMa92},
then $\varphi$ must have the form $\varphi=\frac{1}{3}\chi_{S}$ where $S$ is a
compact subset $\subset\left[  0,3\right]  $ with non-empty interior. It is
determined by the identity%
\begin{equation}
2S=S\cup\left(  S+3\right)  
\label{eqReptileIdentity}%
\end{equation}
(see Figure \ref{ReptileGenealogy}). It follows that $S=\left[  0,3\right]  $.
In fact, iteration of (\ref{eqReptileIdentity}) leads to the following
representation which characterizes points $x$ in $S$:
$x=\sum_{k=1}^{\infty}d_{k}/2^{k}$,
$d_{k}=3\varepsilon_{k}$, $\varepsilon_{k}\in\left\{  0,1\right\}  $.
Hence, using base $2$ for the unit interval $\left[  0,1\right]  $, we get
$S=\left[  0,3\right]  $. The derivation of $\varphi^{\left(  \frac{\pi}%
{2}\right)  }$ from the first Haar wavelet $\varphi^{\left(  \pi\right)
}=\chi_{\left[  0,1\right]  }$ is a special case of the substitution
\begin{equation}
m_{0}\left(  z\right)  \longmapsto m_{0}\left(  z^{3}\right)  ,
\text{\quad or generally,\quad }
m_{0}\left(  z\right)  \longmapsto m_{0}\left(  z^{2p+1}\right)  .
\label{Substitutionz2podd}%
\end{equation}
If $m_{0}$ is an arbitrary low-pass filter with scaling function $\varphi$,
then the argument from Remark \ref{RemMallat} shows that $\widehat
{\varphi_{2p+1}}\left(  \omega\right)  :=\hat{\varphi}\left(  \left(
2p+1\right)  \omega\right)  $ will determine the scaling function for the
substitution $m_{0}\left(  z^{2p+1}\right)  $. Hence%
\begin{equation}
\varphi_{2p+1}\left(  x\right)  =\frac{1}{2p+1}\varphi\left(  \frac{x}%
{2p+1}\right)  ,
\text{\quad and\quad }
\left\|  \varphi_{2p+1}\right\|  _{L^{2}\left(  \mathbb{R}\right)  }=\frac
{1}{\sqrt{2p+1}}\left\|  \varphi\right\|  _{L^{2}\left(  \mathbb{R}\right)  }.
\label{eqHaarSubstitution}%
\end{equation}
In our circular family, we have
$m_{0}^{\left(  \frac{\pi}{2}\right)  }\left(  z\right)
=m_{0}^{\left(  \pi\right)  }\left(  z^{3}\right)  $.
See further discussion in Section \ref{Sing}
and Remark \ref{RemGeometricallyClear}.

\subsubsubsection{\label{Intertwining}Intertwining of the cases $\theta
=\frac{\pi}{2}$ and $\theta=\frac{3\pi}{2}$}

Let us summarize the description of these two representations. By
(\ref{eqPolynomial.9})--(\ref{eqPolynomial.10}) we have%
\begin{equation}
S_{0}^{\frac{3\pi}{2}}e_{n}^{{}} =\frac{1}{\sqrt{2}}
\left( e_{1+2n}+e_{2+2n}\right) ,\qquad S_{1}^{\frac{3\pi}{2}}e_{n}^{{}}
=\frac{1}{\sqrt{2}}\left( -e_{1+2n}+e_{2+2n}\right) ,
\label{eqnu=1.55}%
\end{equation}
and the irreducible invariant subspaces are%
\begin{equation}
\mathcal{H}_{+}^{\frac{3\pi}{2}} =\overline{\operatorname
*{span}}\left\{ e_{-1},e_{0},e_{1},\dots\right\} ,\qquad \mathcal{H}_{-}%
^{\frac{3\pi}{2}} =\overline{\operatorname*{span}}\left\{ e_{-2},e_{-3}%
,\dots\right\} .
\label{eqnu=1.56}%
\end{equation}
Similarly%
\begin{equation}
S_{0}^{\frac{\pi}{2}}e_{n}^{{}} =
\left( e_{2n}+e_{3+2n}\right) /\sqrt{2},
\qquad S_{1}^{\frac{\pi}{2}}e_{n}^{{}} =
\left( e_{2n}-e_{3+2n}\right) /\sqrt{2},
\label{eqnu=1.57}%
\end{equation}
and the associated three irreducible invariant subspaces are%
\begin{equation}
\begin{gathered}
\mathcal{H}_{+}^{\frac{\pi}{2}} =\overline{\operatorname
*{span}}\left\{ e_{0},e_{3},e_{6},\dots\right\} ,\qquad \mathcal{H}_{-}^{\frac
{\pi}{2}} =\overline{\operatorname*{span}}\left\{ e_{-3},e_{-6},e_{-9}%
,\dots\right\} ,\\ \mathcal{J}= \overline{\operatorname*{span}}\left
\{ \dots,e_{-4},e_{-2},e_{-1},e_{1},e_{2},e_{4},\dots\right\} . \end{gathered}
\label{eqnu=1.58}%
\end{equation}

We have noted that $\operatorname*{UHF}\nolimits_{2}$ is not weakly dense in
the last representation but it is so in the first four, so the last
representation cannot be equivalent to any of the former four. Also the
representation on $\mathcal{H}_{+}^{\theta}$ is disjoint from that on
$\mathcal{H}_{-}^{\theta}$ by \cite[Theorem 2.7]{BrJo96b}, for $\theta
=\frac{3\pi}{2}$ and for $\theta=\frac{\pi}{2}$. So the remaining possibility
is that the representation on $\mathcal{H}_{\pm}^{\frac{3\pi}{2}}$ is
unitarily equivalent to that on $\mathcal{H}_{\pm}^{\frac{\pi}{2}}$.
Inspection of the expressions for $S_{i}^{\theta}e_{n}^{{}}$ makes it
plausible that the representation on $\mathcal{H}_{+}^{\frac{3\pi}{2}}$ is
equivalent to that on $\mathcal{H}_{-}^{\frac{\pi}{2}}$, and that that on
$\mathcal{H}_{-}^{\frac{3\pi}{2}}$ is equivalent to that on $\mathcal{H}%
_{+}^{\frac{\pi}{2}}$, and indeed, if one defines an isometry $U$ by%
\begin{equation}
Ue_{n}=e_{-3n-6} \label{eqnu=1.59}%
\end{equation}
then $U|_{\mathcal{H}_{+}^{\frac{3\pi}{2}}}$
from
$\mathcal{H}_{+}^{\frac{3\pi}{2}}$ to $\mathcal{H}_{-}^{\frac{\pi}{2}}$
and
$U|_{\mathcal{H}_{-}^{\frac{3\pi}{2}}}$
from $\mathcal{H}%
_{-}^{\frac{3\pi}{2}}$ to $\mathcal{H}_{+}^{\frac{\pi}{2}}$
are unitary operators, and one computes%
\begin{align}
US_{0}^{\frac{3\pi}{2}}e_{n}^{{}}  &  =\left(
e_{-6n-9}^{{}}+e_{-6n-12}^{{}}\right) /\sqrt{2} 
&  &  =S_{0}^{\frac{\pi}{2}}Ue_{n}^{{}%
},\label{eqnu=1.60}\\
US_{1}^{\frac{3\pi}{2}}e_{n}^{{}}  &  =\left(
-e_{-6n-9}^{{}}+e_{-6n-12}^{{}}\right) /\sqrt{2} 
&  &  =S_{1}^{\frac{\pi}{2}}Ue_{n}%
^{{}}. \label{eqnu=1.61}%
\end{align}
Hence $U$ intertwines the two representations, and if $U$ is restricted to
$\mathcal{H}_{\pm}^{\frac{3\pi}{2}}$ one obtains the expected unitary
intertwiners%
\begin{equation}
U_{1}^{{}}\colon\mathcal{H}_{+}^{\frac{3\pi}{2}}
\longrightarrow\mathcal{H}_{-}^{\frac{\pi}{2}},\qquad
U_{2}^{{}}\colon\mathcal{H}_{-}^{\frac{3\pi}{2}}
\longrightarrow\mathcal{H}_{+}^{\frac{\pi}{2}}.
\label{eqWhyNotIndeed}%
\end{equation}
Now
$U_{1}e_{-1}=e_{-3}$, $U_{2}e_{-2}=e_{0}$,
and hence
\begin{equation}
P\left(  xU_{1}+yU_{2}\right)  P=%
\begin{pmatrix}
0 & 0 & y & 0\\
0 & 0 & 0 & 0\\
0 & 0 & 0 & 0\\
0 & x & 0 & 0
\end{pmatrix}
\label{eqnu=1.64}%
\end{equation}
for $x,y\in\mathbb{C}$, and this is exactly the fixed point set for the map%
\begin{equation}
M_{4}^{{}}\ni A\longmapsto\sum_{i=0}^{1}V_{i}^{\frac{\pi}{2}}AV_{i}%
^{\frac{3\pi}{2}\,\ast}, \label{eqnu=1.65}%
\end{equation}
as it should be by (\ref{eqFinitely.16})--(\ref{eqFinitely.18}). By
\cite[Theorem 5.1]{BJKW97}, these solutions $A$ correspond to operators which
intertwine the two representations.

\subsubsubsection{\label{Sing}Additional remarks on singular points and cycles}

Recall from \cite[Theorem 6.3.6]{Dau92} and \cite[Theorem 3.3.6]{Hor95} that
in order that $\psi_{j,k}\left(  x\right)  =2^{-\frac{j}{2}}\psi\left(
2^{-j}x-k\right)  $ shall be an orthonormal basis for $L^{2}\left(
\mathbb{R}\right)  $ and not merely a tight frame, it is necessary and
sufficient that the set%
\begin{equation}
\left\{  z\in\mathbb{T}\bigm|\left|  m_{0}\left(  z\right)  \right|  =\sqrt
{2}\right\}  =\left\{  z\in\mathbb{T}\mid m_{0}\left(  -z\right)  =0\right\}
\label{eqnu=1.66}%
\end{equation}
does not contain a nontrivial cycle for the doubling map $z\mapsto z^{2}$,
i.e., a finite cyclic subset unequal to $\left\{  1\right\}  $ invariant under
the map $z\mapsto z^{2}$. Inspection of the polynomial $m_{0}^{\left(
\theta\right)  }\left(  z\right)  $ in (\ref{eqnu=1.10}) in the present case
(\ref{eqnu=1.27}) reveals that the condition above is fulfilled for all
$\theta\in\mathbb{T}$ with the sole exception%
\begin{equation}
\theta=\pi/2, \label{eqnu=1.67}%
\end{equation}
where $m_{0}^{\left(  \theta\right)  }\left(  z\right)  $ is given by
(\ref{eqnu=1.47}), and thus the set (\ref{eqnu=1.66}) consists of the three
cube roots of $1$. Indeed, the presence of a nontrivial cycle on $\mathbb{T}$
under $z\mapsto z^{2}$ would imply, by (\ref{eqnu=1.66}) and the fact that
$m_{0}^{\left(  \theta\right)  }$ is a third-degree polynomial, that
$m_{0}^{\left(  \theta\right)  }\left(  \,\cdot\,\right)  $ is a scalar
multiple of $\left(  z+1\right)  \left(  z+e^{i\frac{2\pi}{3}}\right)  \left(
z+e^{-i\frac{2\pi}{3}}\right)  =z^{3}+1$, and this is precisely the case
$\theta=\frac{\pi}{2}$ in (\ref{eqnu=1.67}). (See details of the argument
below.) It is interesting that this is the case of Section
\ref{Casetheta=halfpi} where the decomposition theory of the associated
representation is most singular, being the only case where the restriction of
one of the subrepresentations to $\operatorname*{UHF}\nolimits_{2}$ has a
nontrivial cyclic structure.

Let us give a more detailed justification of the statement above. First note
that cycles on $\mathbb{T}$ are not subgroups of $\mathbb{T}$ but rather
cyclic orbits on $\mathbb{T}$ under the $z\mapsto z^{2}$ action of one of the
cyclic groups $\mathbb{Z}_{k}$, $k=1,2,\dots$. Such a cyclic orbit $C_{k}$
with $k$ distinct points $z_{1},\dots,z_{k}$ must be of the form
$z_{1}\rightarrow z_{2}\rightarrow\dots\rightarrow z_{k}\rightarrow z_{1}$,
where $z_{i+1}^{{}}=z_{i}^{2}$ if $i=1,2,\dots,k-1$, and $z_{k}^{2}=z_{1}^{{}%
}$. Hence points $c$ in an orbit $C_{k}$ must satisfy $c^{2^{k}}=c$, and each
$c$ must be a $\left(  2^{k}-1\right)  $'th root of $1$. Different orbits must
be disjoint, and their union will be invariant under $z\mapsto z^{2}$ acting
on $\mathbb{T}$. The converse is not true. For example, the subset $\left\{
1,-1\right\}  \subset\mathbb{T}$ is invariant under $z\mapsto z^{2}$ while not
a cycle, and not even the union of cycles. Note also that we can have
different $\left(  2^{k}-1\right)  $'th roots $c$ of $1$ defining different
cyclic orbits for the same $k$. If $k=1$ or $k=2$, then in each case there is
only one orbit, but if $k=3$, there are two choices. Since $m_{0}^{\left(
\theta\right)  }$ for each $\theta$ is a polynomial of degree $3$, the
cardinality of a cycle contained in the set (\ref{eqnu=1.66}) is at most $3$.
Thus, if $z$ is contained in such a cycle, we must have one of the
possibilities $z^{2}=z$, $z^{4}=z$, $z^{8}=z$. Hence the cycles of length at
most $3$ are $\left\{  1\right\}  $, $\left\{  \omega,\omega^{2}\right\}  $
where $\omega=e^{i\frac{2\pi}{3}}$, $\left\{  \rho,\rho^{2},\rho^{4}\right\}
$ where $\rho=e^{i\frac{2\pi}{7}}$, and $\left\{  \bar{\rho},\bar{\rho}%
^{2},\bar{\rho}^{4}\right\}  =\left\{  \rho^{6},\rho^{5},\rho^{3}\right\}  $.
But as $m_{0}\left(  -1\right)  =0$ always, $\left(  z+1\right)  $ is always a
factor of $m_{0}\left(  z\right)  $, and since the cycle should be different
from the trivial cycle $\left\{  1\right\}  $, we are reduced to the case
$\left\{  \omega,\omega^{2}\right\}  $. The other cycles would make
$m_{0}^{\left(  \theta\right)  }$ divisible by a polynomial of degree at least
$4$, which of course is impossible. Thus we are left with the case%
\begin{equation}
m_{0}\left(  z\right)  =\frac{1}{\sqrt{2}}\prod_{k=0}^{2}\left(  \omega
^{k}+z\right)  =\frac{1}{\sqrt{2}}\left(  1+z^{3}\right)  , \label{eqomegod}%
\end{equation}
which is exactly the case $\theta=\frac{\pi}{2}$.

Note, more generally, that the wavelets which arise from substitutions, as
defined in (\ref{Substitutionz2podd})--(\ref{eqHaarSubstitution}) with filter
function $m_{0}^{\left(  p\right)  }\left(  z\right)  =m_{0}^{{}}\left(
z_{{}}^{2p+1}\right)  $, will have those additional cycles $C_{k}$ which are
contained in the $\left(  2p+1\right)  $'th roots of $1$, $\left\{
z\in\mathbb{T}\mid z^{2p+1}=1\right\}  $.
We will show in a forthcoming paper that
this leads to a decomposition of the
representation of $\mathcal{O}_{2}$ associated to $m_0^{(p)}$
over the new cycles.

It is interesting to note that the same cycles as described above arise in a
different context in \cite{BrJo96b} in connection with a family of discrete
series of representations of $\mathcal{O}_{N}$. These representations are
called permutative representations, and the cycles represent the finite
decompositions of irreducible representations of $\mathcal{O}_{N}$ when
restricted to $\operatorname*{UHF}\nolimits_{N}$.%

\subsubsubsection{\label{Symmetry}The symmetry $\theta\protect\mapsto
\pi-\theta$}

Note that the two points $\theta=0$ and $\theta=\pi$ are interesting in that
the representation theory is regular, but these points correspond to mother
and father functions which are simple rescalings of those of the Haar wavelet
(see Figure \ref{FatherMothertheta=khalfpi}):
\begin{equation}
\theta=0 \colon
\begin{cases}
m_{0}\left(  z\right)   &  =\left(  z^{2}+z^{3}\right) /\sqrt{2} \\
m_{1}\left(  z\right)   &  =\left(  1-z\right) /\sqrt{2}
\end{cases}
\qquad
\theta=\pi \colon
\begin{cases}
m_{0}\left(  z\right)   &  =\left(  1+z\right) /\sqrt{2} \\
m_{1}\left(  z\right)   &  =\left(  z^{2}-z^{3}\right) /\sqrt{2}
\end{cases}
\label{eqnu=1.68and69}%
\end{equation}
Thus the representation of $\mathcal{O}_{2}$ is very sensitive to simple
rescaling of $\varphi$ and $\psi$. In fact the mother function $\psi$ is the
same in the two cases $\theta=0$ and $\theta=\pi$, and this common $\psi$ has
the following symmetry property $\psi\left(  3-x\right)  =-\psi\left(
x\right)  $, which in turn is a special case of a more general reflection
symmetry (\ref{eqMirrorMother}) to be discussed in Proposition \ref{ProMirror}%
(\ref{ProMirror(a)}) below.

The symmetry $0\mapsto\pi$ is a special case of a symmetry $\theta\mapsto
\pi-\theta$, which we will now analyze further. If this transformation is
substituted in (\ref{eqnu=1.27}), we note that it corresponds to the following
reversal:%
\begin{equation}
\left(  a_{0},a_{2},a_{2},a_{3}\right)  \longmapsto\left(  a_{3},a_{2}%
,a_{1},a_{0}\right)  ; \label{eqSwapCoefficients}%
\end{equation}
or equivalently,
\begin{equation}
m_{0}^{\left(  \pi-\theta\right)  }\left(  z\right)  =m_{0}^{\left(
\theta\right)  }\left(  z^{-1}\right)  z^{3}. \label{eqReversal}%
\end{equation}
The following proposition shows that the $\theta\mapsto\pi-\theta$ reflection
applied to $m_{0}^{\left(  \theta\right)  }$ implements the $x\mapsto3-x$
transformation on the scaling function $\varphi$
(see Figure \ref{FatherMotherNear}a,b).
It is interesting to note that, despite this left-right
mirror symmetry of the graphs in the family of scaling functions
$\varphi^{\left(  \theta\right)  }$, the two associated representations of
$\mathcal{O}_{2}$ on $L^{2}\left(  \mathbb{T}\right)  $ which correspond,
respectively, to $\theta$ and $\pi-\theta$, are \emph{not} unitarily
equivalent, by the results above, except of course at the two fixed points
$\frac{\pi}{2}$ and $\frac{3\pi}{2}$ for $\theta\mapsto\pi-\theta$, where the
representation theory also happens to be exceptional. See subsections
\ref{Casetheta=3halfpi} and \ref{Casetheta=halfpi} above.

\begin{proposition}
\label{ProMirror}\renewcommand{\theenumi}{\alph{enumi}}Let $m_{0}^{\left(
\theta\right)  }\left(  z\right)  $ be the filter functions indexed by
$\theta$ and corresponding to the given coefficients in the family
\textup{(\ref{eqnu=1.27})}. Let $\varphi^{\left(  \theta\right)  }\left(
x\right)  $ be the associated scaling function \textup{(}alias, father
function\textup{)} and $\psi^{\left(  \theta\right)  }$ the mother function
corresponding to the pair $\left(  m_{0}^{\left(  \theta\right)  }%
,m_{1}^{\left(  \theta\right)  }\right)  $ of low/high-pass wavelet filters.
Let $S_{i}^{\left(  \theta\right)  }$ be the corresponding operators from
\textup{(\ref{eqPolynomial.4}\textup{)}}.

\begin{enumerate}
\item \label{ProMirror(a)}The symmetry relations%
\begin{align}
\varphi^{\left(  \pi-\theta\right)  }\left(  x\right)   &  =\varphi^{\left(
\theta\right)  }\left(  3-x\right)  ,\quad\theta\in\left[  -\pi,\pi\right]
,\;x\in\mathbb{R},\label{eqMirrorFather}\\
\psi^{\left(  \pi-\theta\right)  }\left(  x\right)   &  =-\psi^{\left(
\theta\right)  }\left(  3-x\right)  , \label{eqMirrorMother}%
\end{align}
are valid.\smallskip

\item \label{ProMirror(b)}The corresponding representations $\bp^{\left(
\theta\right)  }$ and $\bp^{\left(  \pi-\theta\right)  }$ \textup{(}given
by $\bp_{{}}^{\left(  \theta\right)  }\left(  s_{i}^{{}}\right)
=S_{i}^{\left(  \theta\right)  }$\textup{)} satisfy%
\begin{equation}
W\bp^{\left(  \theta\right)  }=\smash[t]{\left(  \bp^{\left(  \pi-\theta\right)  }
\circ\tau_{\left(
\begin{smallmatrix}
1 & 0\\
0 & -1
\end{smallmatrix}
\right)  }\right)  }W, \label{eqWinter}%
\end{equation}
where $\left(  Wf\right)  \left(  z\right)  =z^{-3}f\left(  z^{-1}\right)  $,
and $\smash[b]{\tau_{\left(
\begin{smallmatrix}
1 & 0\\
0 & -1
\end{smallmatrix}
\right)  }}$ is the automorphism of $\mathcal{O}_{2}$ given in
\textup{(\ref{eqFinitely.3}\textup{)}} for $g=\left(
\begin{smallmatrix}
1 & 0\\
0 & -1
\end{smallmatrix}
\right)  $.
\end{enumerate}
\end{proposition}

\begin{proof}
Introducing $z=e^{-i\omega}$, $\omega\in\mathbb{R}$, the identity
(\ref{eqReversal}) above reads
\begin{equation}
m_{0}^{\left(  \pi-\theta\right)  }\left(  \omega\right)  =e^{i3\omega}%
m_{0}^{\left(  \theta\right)  }\left(  -\omega\right)  =e^{i3\omega}%
\overline{m_{0}^{\left(  \theta\right)  }\left(  \omega\right)  },\quad
\omega\in\mathbb{R}. \label{eqMirroromega}%
\end{equation}
Generally for third degree, the correspondence $m_{0}\leftrightarrow\varphi$
is given by the following functional identity in $L^{2}\left(  \mathbb{R}%
\right)  $:%
\begin{equation}
\varphi\left(  x/2\right) /\sqrt{2} =a_{0}\varphi\left(
x\right)  +a_{1}\varphi\left(  x-1\right)  +a_{2}\varphi\left(  x-2\right)
+a_{3}\varphi\left(  x-3\right)  , \label{eqMirrorFunctional}%
\end{equation}
and the boundary conditions, $\varphi\left(  0\right)  =\varphi\left(
3\right)  =0$, i.e., $\varphi$ is uniquely determined by these conditions and
the normalization $\hat{\varphi}\left(  0\right)  =\left(  2\pi\right)
^{-\frac{1}{2}}$. See \cite{Pol92} for details. This applies to both the pair
$\left(  m_{0}^{\left(  \theta\right)  },\varphi^{\left(  \theta\right)
}\right)  $ and the pair $\left(  m_{0}^{\left(  \pi-\theta\right)  }%
,\varphi^{\left(  \pi-\theta\right)  }\right)  $, so we get%
\begin{multline}
\varphi^{\left(  \pi-\theta\right)  }\left(  x/
2\right) /\sqrt{2}\label{eqMirrorpiminustheta}\\
=a_{3}\varphi^{\left(  \pi-\theta\right)  }\left(  x\right)  +a_{2}%
\varphi^{\left(  \pi-\theta\right)  }\left(  x-1\right)  +a_{1}\varphi
^{\left(  \pi-\theta\right)  }\left(  x-2\right)  +a_{0}\varphi^{\left(
\pi-\theta\right)  }\left(  x-3\right)  .
\end{multline}
As noted, $\varphi^{\left(  \pi-\theta\right)  }\left(  \,\cdot\,\right)  $ is
the unique normalized $L^{2}\left(  \mathbb{R}\right)  $-solution to this
identity, subject to $\varphi^{\left(  \pi-\theta\right)  }\left(  0\right)
=\varphi^{\left(  \pi-\theta\right)  }\left(  3\right)  =0$. But, if
$\varphi^{\left(  \theta\right)  }$ is the solution corresponding to
$m_{0}^{\left(  \theta\right)  }$, then a direct substitution $x\mapsto6-x$
shows that the mirrored function $x\mapsto\varphi^{\left(  \theta\right)
}\left(  3-x\right)  $ satisfies (\ref{eqMirrorpiminustheta}), and we conclude
from the uniqueness that%
\begin{equation}
\varphi^{\left(  \pi-\theta\right)  }\left(  x\right)  =\varphi^{\left(
\theta\right)  }\left(  3-x\right)  ,\quad x\in\mathbb{R},
\label{eqMirrorConclusion}%
\end{equation}
as claimed in the Proposition. The proof of (\ref{eqMirrorMother}) is similar,
or see Remark \ref{RemMallat} below. We resume the proof of Proposition
\ref{ProMirror}(\ref{ProMirror(b)}) after the following remark.%
\renewcommand{\qed}{}%
\end{proof}

\begin{remark}
\label{RemMallat}Proposition \textup{\ref{ProMirror}%
\textup{(\ref{ProMirror(a)})}} may alternatively be proved from the Mallat
algorithm as follows: If $\varphi^{\left(  \theta\right)  }$, $\psi^{\left(
\theta\right)  }$ are the father and mother functions at the angle $\theta$,
and the transformation $\theta\mapsto\pi-\theta$ is used on
\textup{(\ref{eqnu=1.27})}, we obtain \textup{(\ref{eqSwapCoefficients}%
\textup{)}} and \textup{\textup{(\ref{eqReversal})}} as before, i.e.,%
\begin{equation}
m_{0}\left(  z\right)   \longmapsto z^{3}\overline{m_{0}\left(  z\right)
}=m_{1}\left(  -z\right) 
\text{\quad and\quad }
m_{1}\left(  z\right)   \longmapsto m_{0}\left(  -z\right)  .
\label{eqSwapm01and10}%
\end{equation}
Applying the Mallat algorithm
$\hat{\varphi}\left(  t\right)  =\left(  2\pi\right)  ^{-\frac{1}{2}}
\prod_{k=1}^{\infty}\left(  m_{0}\left(  e^{-it2^{-k}
}\right)  /\sqrt{2}\right)  $,
we obtain%
\begin{equation}
\widehat{\varphi^{\left(  \pi-\theta\right)  }}\left(  t\right)
=e^{-i3t}\widehat{\varphi^{\left(  \theta\right)  }}\left(  -t\right)  ,
\label{eqMallatFatherhat}%
\end{equation}
and thus by Fourier transform,%
\begin{equation}
\varphi^{\left(  \pi-\theta\right)  }\left(  x\right)  =\varphi^{\left(
\theta\right)  }\left(  3-x\right)  , \label{eqMirrorFourierFather}%
\end{equation}
which is \textup{(\ref{eqMirrorFather})}. On the other hand,%
\begin{equation}
\psi\left(  x\right)   =\sqrt{2}\sum_{k}\left(  -1\right)  ^{k}%
a_{3-k}\varphi\left(  2x-k\right)  ,
\label{eqFatherbyMother}%
\end{equation}
and so
\begin{multline}
\psi^{\left(  \pi-\theta\right)  }\left(  x\right)   =\sqrt{2}\sum
_{k}\left(  -1\right)  ^{k}a_{3-k}^{\left(  \pi-\theta\right)  }%
\varphi^{\left(  \pi-\theta\right)  }\left(  2x-k\right) 
\label{eqMirrorIndeedMotherConclusion}\\
=\sqrt{2}\sum_{k}\left(  -1\right)  ^{k}a_{k}^{\left(  \theta\right)
}\varphi^{\left(  \theta\right)  }\left(  3-\left(  2x-k\right)  \right)  \\
=\sqrt{2}\sum_{k}\left(  -1\right)  ^{3-k}a_{3-k}^{\left(  \theta\right)
}\varphi^{\left(  \theta\right)  }\left(  2\left(  3-x\right)  -k\right)
=-\psi^{\left(  \theta\right)  }\left(  3-x\right)  ,
\end{multline}
which is \textup{(\ref{eqMirrorMother})}.

It is important to note that the infinite product argument
works
even if
$\varphi^{\left(  \theta\right)  }\left(  x\right)  $ is not continuous in
$x$. Since
\begin{equation}
\left|  m_{0}^{\left(  \theta\right)  }\left(  \omega\right)  \right|
^{2}+\left|  m_{0}^{\left(  \theta\right)  }\left(  \omega+\pi\right)
\right|  ^{2}=2,\quad\omega\in\mathbb{R}, \label{eqSumSquarem0}%
\end{equation}
it is known that 
the infinite products%
\begin{equation}
\left(  2\pi\right)  ^{-\frac{1}{2}}\prod_{k=1}^{\infty}2^{-\frac{1}{2}}%
m_{0}^{\left(  \theta\right)  }\left(  \frac{\omega}{2^{k\mathstrut}}\right)
,\quad\left(  2\pi\right)  ^{-\frac{1}{2}}2^{-\frac{1}{2}}m_{1}^{\left(
\theta\right)  }\left(  \frac{\omega}{2}\right)  \prod_{k=2}^{\infty}%
2^{-\frac{1}{2}}m_{0}^{\left(  \theta\right)  }\left(  \frac{\omega
}{2^{k\mathstrut}}\right)  \label{eqFatherhatomega}%
\end{equation}
are well defined and represent $\widehat{\varphi^{\left(
\theta\right)  }}$, $\widehat{\psi^{\left(  \theta\right)  }}$, where
$\varphi^{\left(  \theta\right)  },\psi^{\left(  \theta\right)  }\in
L^{2}\left(  \mathbb{R}\right)  $ \cite{Dau92}.
\end{remark}

\begin{proof}
[Proof of Proposition \textup{\ref{ProMirror}(\ref{ProMirror(b)})}]Let us
consider the two operators $S_{0}^{\left(  \theta\right)  }$ and
$S_{0}^{\left(  \pi-\theta\right)  }$ in $L^{2}\left(  \mathbb{T}\right)  $
individually, and as part of a pair of $\mathcal{O}_{2}$-representations.
While the two $\mathcal{O}_{2}$-representations are inequivalent, the two
$S_{0}$-operators alone are unitarily equivalent. This follows from the
general fact that any operator of the form (\ref{eqIntro.12}) coming from a
wavelet is unitarily equivalent to the shift of infinite multiplicity by
\cite[Lemma 9.3]{BrJo97b}. The explicit intertwiner can also be calculated
directly as follows: Let $m\left(  z\right)  =a_{0}+a_{1}z+\dots+a_{D}z^{D}$,
$m^{\prime}\left(  z\right)  :=z^{D}m\left(  z^{-1}\right)  $, and define
three operators $S$, $S^{\prime}$, and $W$
(acting on $f\in L^{2}\left(  \mathbb{T}\right)  $)
by%
\[
Sf\left(  z\right) :=m\left(  z\right)  f\left(  z^{2}\right)  ,\quad
S^{\prime}f\left(  z\right) :=m^{\prime}\left(  z\right)  f\left(
z^{2}\right)  ,\quad
Wf\left(  z\right) :=z^{-D}f\left(  z^{-1}\right)  .
\]
Then $W\colon L^{2}\left(  \mathbb{T}\right)  \rightarrow L^{2}\left(
\mathbb{T}\right)  $ is a unitary intertwining operator for
$S$ and $S^{\prime}$, i.e.,%
\[
WS=S^{\prime}W
\]
holds, as can be verified by a direct calculation. The minus sign in the
second symmetry formula (\ref{eqMirrorMother}) is still reflected in the
$\mathcal{O}_{2}$-representations as follows. Let $D=3$ and $m=m_{0}^{\left(
\theta\right)  }$, and consider the two $\mathcal{O}_{2}$-representations
$\bp^{\left(  \theta\right)  }$, $\bp^{\left(  \pi-\theta\right)  }$, $i=0,1$.
We then have $m_{1}^{\left(  \pi-\theta\right)  }\left(  z\right)  =-z_{{}%
}^{3}m_{1}^{\left(  \theta\right)  }\left(  z_{{}}^{-1}\right)  $, and thus%
\begin{equation}
WS_{0}^{\left(  \theta\right)  }  =S_{0}^{\left(  \pi-\theta\right)
}W,\qquad
WS_{1}^{\left(  \theta\right)  }  =-S_{1}^{\left(  \pi-\theta\right)  }W,
\label{eqWinterplusminus}%
\end{equation}
where again $S_{i}^{\left(  \theta\right)  }=\bp_{{}}^{\left(  \theta\right)
}\left(  s_{i}^{{}}\right)  $. Hence $W$ intertwines the $\theta
$-representation $\bp^{\left(  \theta\right)  }$ with the $\left(  \pi
-\theta\right)  $-representation $\bp^{\left(  \pi-\theta\right)  }$, modified
by the automorphism of $\mathcal{O}_{2}$ induced by $g=\left(
\begin{smallmatrix}
1 & 0\\
0 & -1
\end{smallmatrix}
\right)  \in\mathrm{U}\left(  2\right)  $; see (\ref{eqFinitely.3}).
\end{proof}

\subsubsubsection{\label{Continuity}Continuity of scaling functions: Stability interval}

Historically the special case $\theta=\frac{7\pi}{6}$ in (\ref{eqnu=1.27}) was
discovered first. In that case,
\[
a_{0}=\frac{1+\sqrt{3}}{4\sqrt{2}},\quad a_{1}=\frac{3+\sqrt{3}}{4\sqrt{2}%
},\quad a_{2}=\frac{3-\sqrt{3}}{4\sqrt{2}},\quad a_{3}=\frac{1-\sqrt{3}%
}{4\sqrt{2}},
\]
which is the (by now) well known Daubechies wavelet; see Figure
\ref{FatherMotherNear}a,b \cite[Chapter 6]{Dau92}. It was analyzed further in
\cite{Pol92}, where it was shown to have scaling function $\varphi\left(
\,\cdot\,\right)  $ continuous and one-sided differentiable in $x$, support on
$\left[  0,3\right]  $, $\varphi\left(  0\right)  =\varphi\left(  3\right)
=0$. It is left-differentiable at every dyadic $x$, but it is not
right-differentiable at any dyadic $x$ in $\left[  0,3\right\rangle $. Since
Daubechies established continuity by a matrix spectral estimate (see
\cite[Theorem 7.2.1]{Dau92}), it follows from her estimates that the scaling
function $\varphi^{\left(  \theta\right)  }\left(  x\right)  $ will also be
continuous in an open interval containing $\theta=\frac{7\pi}{6}$. It is
interesting to note that Daubechies's spectral estimation involves the two
matrices $V_{i}^{\ast}$, $i=0,1$, given in (\ref{eqnu=1.28and29}%
) above. The discussion in our previous section indicates that the stability
interval in the $\theta$ variable, $\theta_{0}<\theta<\theta_{1}$, must have
$\pi<\theta_{0}$ and $\theta_{1}<\frac{3\pi}{2}$.

The pictures of the scaling function in this paper are generated with the aid
of the cascade algorithm described in \cite[Section 6.5]{Dau92}.

For uniform convergence of the cascade approximants to $\varphi$, one has to
assume that $\varphi$ is H\"older continuous \cite[Proposition 6.5.2]{Dau92}.
Figure \ref{ReptileGenealogy} shows clearly that this uniform convergence may
fail abysmally even when $\varphi$ is a simple step function. However, we see
from Figure \ref{ReptileGenealogy} that the cascade approximants converge in
the distribution sense, and even in the weak-$L^{2}$ sense, to $\varphi$ when
$\theta=\frac{\pi}{2}$.

The other assumption in
Daubechies's cascade approximation
is the orthogonality of $\mathbb{Z}$-translates, in the
form \cite[(6.5.4)--(6.5.5), p.\ 204]{Dau92},
and,
as we will discuss in Remark \ref{RemGeometricallyClear} below,
that fails when $\theta=\frac{\pi}{2}$, but
is satisfied at all other
values of $\theta$ by Section \ref{Sing} above.

More importantly, Daubechies
states
in \cite[Chapter 6 footnote 9 and Section 6.3]{Dau92} that,
even if
$\varphi$ is not assumed continuous,
we still have $L^{2}\left(\mathbb{R}\right)$ norm convergence
of the cascade-algorithm approximation,
as long as the $\mathbb{Z}$-translates
are mutually orthogonal; and, as
we noted, this
orthogonality holds
whenever
$\theta\neq\frac{\pi}{2}$.
This will be discussed in a forthcoming paper \cite{BrJo98b}.

\begin{remark}
\label{RemGeometricallyClear}
Since
$\varphi^{\left(\frac{\pi}{2}\right)}=\frac{1}{3}\chi_{\left[0,3\right]}$,
it is geometrically clear that the $\mathbb{Z}$-translates of 
$\varphi^{\left(\frac{\pi}{2}\right)}$ in
$L^{2}\left(\mathbb{R}\right)$ will not be
mutually orthogonal \textup{(}see Figure \textup{\ref{ReptileGenealogy})},
and we have shown in Section \textup{\ref{Sing}}
that $\varphi^{\left(\frac{\pi}{2}\right)}$ is the
unique scaling function
in the family
$\left\{\varphi^{\left(\theta\right)}\right\}$
which does not have orthogonal
$\mathbb{Z}$-translates.
The cascade
algorithm, which is used in generating the present graphics, is based on
an iteration of
\textup{(\ref{eqScalingRelation1})}
but is also closely connected to iteration of
$F_{0}^{*}$ in \textup{(\ref{eqQuadmf})}. Let
\begin{equation}
\label{CoefficientAsIntegral}
c_{k}:=\frac{1}{\sqrt{2}}\int_{\mathbb{R}}
\overline{\varphi\left(x-k\right)}\varphi\left(\frac{x}{2}\right)\,dx.
\end{equation}
In the case when
$\left\{\varphi\left(\,\cdot\,-k\right)\right\}_{k\in\mathbb{Z}}$
is an orthonormal basis,
we get $c_{k}=a_{k}$, $k\in\mathbb{Z}$,
by \textup{(\ref{eqScalingRelation1})}; but, in
general, we have a
discrepancy $c_{k}\neq a_{k}$ which
leads to a rather poor
approximation with $a_{k}$-cascades.
For a more explicit estimate we need the following:
\begin{lemma}
\label{LemDiscrepancyEstimate}
Let $m_{0}$ be a
low-pass wavelet filter with
corresponding scaling function $\varphi$
and suppose
that
the $\mathbb{Z}$-translates of
$\varphi$ are orthogonal. Let $\varphi_{p}$
be the scaling function corresponding
to the substitution
$m_{0}\left(z^{2p+1}\right)$,
and let
\[c_{k}^{\left(p\right)}:=\frac{1}{\sqrt{2}}\int_{\mathbb{R}}
\overline{\varphi_{p}\left(x-k\right)}\varphi_{p}\left(\frac{x}{2}\right)\,dx.
\]
Then
\begin{equation}
\label{eqDiscrepancyEstimate}
\sum_{k}\left|c_{k}^{\left(p\right)}\right|^{2}\leq\frac{1}{2p+1}.
\end{equation}
\end{lemma}
\begin{proof}
{}From \cite{Dau92} or \cite[Proposition 12.4]{BrJo97b},
we have
\begin{equation}
\label{eqPhiHatSum}
\sum_{l\in\mathbb{Z}}\left|\hat{\varphi}\left(\omega+2\pi l\right)\right|^{2}
\equiv\frac{1}{2\pi}.
\end{equation}
Since $\widehat{\varphi_{p}}\left(\omega\right)
=\hat{\varphi}\left(\left(2p+1\right)\omega\right)$,
we conclude that
\begin{equation}
\label{eqPhiHatBound}
\sum_{l}\left|\hat{\varphi}\left(
\left(2p+1\right)
\left(\omega+2\pi l\right)
\right)\right|^{2}
\leq\frac{1}{2\pi}.
\end{equation}
This second summation is just
one of the $2p+1$ residue classes
for the full $\mathbb{Z}$ summation in \textup{(\ref{eqPhiHatSum})}.
But the formula for $c_{k}$ yields
\begin{multline*}
\sum_{k\in\mathbb{Z}}\left|c_{k}^{\left(p\right)}\right|^{2}
=\sum_{k\in\mathbb{Z}}\left|
\int_{0}^{2\pi}e^{ik\omega}m_{0}\left(
\left(2p+1\right)
\omega\right)
\sum_{l\in\mathbb{Z}}\left|\hat{\varphi}\left(
\left(2p+1\right)
\left(\omega+2\pi l\right)
\right)\right|^{2}\,d\omega
\right|^{2}   \\
=2\pi\int_{0}^{2\pi}\left|m_{0}\left(
\left(2p+1\right)
\omega\right)\right|^{2}
\left(\sum_{l}\left|\hat{\varphi}\left(
\left(2p+1\right)
\left(\omega+2\pi l\right)
\right)\right|^{2}
\right)^{2}
\,d\omega   \\
=\frac{1}{\left(2p+1\right) 2\pi}\int_{0}^{2\pi}\left|m_{0}\left(
\omega\right)\right|^{2}
\sum_{j=0}^{2p}
\left(\sum_{l\in\mathbb{Z}}2\pi\left|\hat{\varphi}\left(
\omega+\left(j+\left(2p+1\right)
l\right)
2\pi \right)\right|^{2}
\right)^{2}
\,d\omega   \\
\leq\frac{1}{\left(2p+1\right) 2\pi}\int_{0}^{2\pi}\left|m_{0}\left(
\omega\right)\right|^{2}
\sum_{j=0}^{2p}
\sum_{l\in\mathbb{Z}}2\pi\left|\hat{\varphi}\left(
\omega+
\smash{\underbrace{\left(j+\left(2p+1\right)l\right)}_{n}}
2\pi \right)\right|^{2}
\vphantom{\underbrace{\left(j+\left(2p+1\right)l\right)}_{n}}
\,d\omega   \\
=\frac{1}{\left(2p+1\right) 2\pi}\int_{0}^{2\pi}\left|m_{0}\left(
\omega\right)\right|^{2}
\underbrace{\sum_{n\in\mathbb{Z}}2\pi\left|\hat{\varphi}\left(
\omega+n\cdot 2\pi \right)\right|^{2}}_{1}
\,d\omega   \\
=\frac{1}{2p+1}\cdot\frac{1}{2\pi}\int_{0}^{2\pi}
\left|m_{0}\left(\omega\right)\right|^{2}\,d\omega  
=\sum_{k\in\mathbb{Z}}\left|a_{k}\right|^{2}\frac{1}{2p+1}=\frac{1}{2p+1}.
\hbox to2em{\hfill\hbox to8pt{\qed}}
\end{multline*}
\renewcommand{\qed}{}
\end{proof}

\TeXButton{Figure}{\begin{figure}
\newcounter{coin}
\setcounter{coin}{-4}
\setlength{\unitlength}{10pt}
\begin{picture}(33,15)(-12,-3)
\put(-12,0){\line(1,0){33}}
\multiput(-9,0)(3,0){10}{\line(0,1){0.3}}
\multiput(-10,-2)(3,0){10}{\addtocounter{coin}{1}\makebox(2,2){\arabic{coin}}}
\put(19,-2){\makebox(2,2){$k$}}
\setcounter{coin}{-4}
\multiput(-9,0)(3,1){3}{\circle*{0.3}}
\multiput(-10,0)(3,1){3}{\makebox(2,2){\addtocounter
{coin}{1}$c_{\arabic{coin}}$}}
\multiput(0,3)(3,0){3}{\circle*{0.3}}
\put(-2,3){\makebox(2,2){\addtocounter
{coin}{1}$c_{\arabic{coin}}$}}
\multiput(2,3)(3,0){2}{\makebox(2,2){\addtocounter
{coin}{1}$c_{\arabic{coin}}$}}
\multiput(9,3)(3,-1){4}{\circle*{0.3}}
\put(9,3){\makebox(2,2){\addtocounter
{coin}{1}$c_{\arabic{coin}}$}}
\multiput(11,2)(3,-1){3}{\makebox(2,2){\addtocounter
{coin}{1}$c_{\arabic{coin}}$}}
\put(-9,0){\line(3,1){9}}
\put(0,3){\line(1,0){9}}
\put(9,3){\line(3,-1){9}}
\multiput(0,9)(9,0){2}{\circle{0.3}}
\put(-1,9){\makebox(2,2){$a_{0}$}}
\put(8,9){\makebox(2,2){\rlap{$a_{3}=\frac{1}{\sqrt{2}}$}\hphantom{$a_{3}$}}}
\setcounter{coin}{0}
\multiput(3,0)(3,0){2}{\circle{0.3}}
\multiput(2,0)(3,0){2}{\makebox(2,2){\addtocounter
{coin}{1}$a_{\arabic{coin}}$}}
\multiput(0,0)(9,0){2}{\line(0,1){9}}
\end{picture}
\renewcommand{\frogcapCorrelationCoefficients}{: The
correlation coefficients $c_{k}$, for which
$\sum_{k=-2}^{5}c_{k}^{2}=\frac{23}{81}<\frac{1}{3}$,
as compared to the scaling coefficients $a_{k}$,
for which $\sum_{k=0}^{3}a_{k}^{2}=1$.}
\caption{Correlation coefficients for
$\theta=\frac{\pi}{2}$\protect\frogcapCorrelationCoefficients}
\label{CorrelationCoefficients}
\end{figure}}

We now illustrate this for $\varphi^{\left(\frac{\pi}{2}\right)}$.
Since $\varphi^{\left(\frac{\pi}{2}\right)}
=\frac{1}{3}\chi_{\left[0,3\right]}$, it is easy
to compute exactly the correlation coefficients
$c_{k}$ of \textup{(\ref{CoefficientAsIntegral})}.
The nonzero coefficients are:
\[
c_{-2} = c_{5} = 1/9\sqrt{2},   \quad
c_{-1} = c_{4} = 2/9\sqrt{2}, \text{\quad and\quad }  
c_{0} = c_{1} = c_{2} = c_{3} = 1/3\sqrt{2},
\]
which should be compared with \textup{(\ref{eqnu=1.46})}. They are
also illustrated in Figure \textup{\ref{CorrelationCoefficients}},
and a comparison with Figure \textup{\ref{ReptileGenealogy}} suggests
that replacing the $a_{k}$'s in the
cascades with the $c_{k}$'s might
possibly lead to a better
approximation. Good approximations
are not known in the non-orthogonal
case. For more details, see \cite[pp.\ 204--206]{Dau92}.
The $c_{k}$ numbers are those
which may be inserted into
the spline approximation that
is also discussed in \cite[pp.\ 206--207]{Dau92}
to build in tight frame parameters
in the approximation.

The problem with this
substitution of the $c_{k}$'s into
the
cascade
algorithm is that, in
the non-orthogonal case, we
will have \textup{(}see Lemma \textup{\ref{LemDiscrepancyEstimate})}
$\sum_{k}\left|c_{k}\right|^{2}<1$. Compare this
to the normalization property
$\sum_{k}\left|a_{k}\right|^{2}=1$ from \textup{(\ref{eqnu=1.4and5})},
or \textup{(\ref{eqnu=1.12and13})} in the special case.
\end{remark}

\section{Conclusions}

We have demonstrated how a representation-theoretic approach
to the construction of compactly supported wavelets
in $\mathbb{R}^{d}$ leads to:
\begin{enumerate}
\item a coordinate-free display of the examples,
\item a finite-dimensional matrix algorithm for computing irreducibility properties,
\item a formula for decomposition into orthogonal sums of irreducibles.
\end{enumerate}
The theory is illustrated in the simplest cases where the power of the
representation-theoretic approach comes into play.

\begin{acknowledgements}
We are indebted to Rune Kleveland for assistance with the computations in
Section \textup{\ref{Casenu=1N=2d=2}}, to Erik L\o w for informing us about
the cascade algorithm used in the graphics, to the University of Oslo for
support, and to Brian Treadway for excellent typesetting, graphics,
discovering the symmetry $\theta\mapsto\pi-\theta$, and manuscript coordination.
\end{acknowledgements}

\bibliographystyle{amsplain}
\bibliography{jorgen}

\providecommand{\bysame}{\leavevmode\hbox to3em{\hrulefill}\thinspace}
\begin{thebibliography}{10}

\bibitem{Ban91}
Christoph Bandt, \emph{Self-similar sets \textup{5:} {I}nteger matrices and
  fractal tilings of $\mathbb{R}^n$}, Proc. Amer. Math. Soc. \textbf{112}
  (1991), 549--562.

\bibitem{Ban97}
\bysame, \emph{Self-similar tilings and patterns described by mappings}, The
  Mathematics of Long-Range Aperiodic Order (Waterloo, Ontario, 1995)
  (R.~Moody, ed.), NATO Adv. Sci. Inst. Ser. C Math. Phys. Sci., vol. 489,
  Kluwer Academic Publishers Group, 1997.

\bibitem{BrJo98b}
O.~Bratteli and P.E.T. Jorgensen, \emph{Convergence of the cascade algorithm at
  irregular scaling functions}, in preparation.

\bibitem{BrJo96b}
\bysame, \emph{Iterated function systems and permutation representations of the
  {C}untz algebra}, Mem.\ Amer.\ Math.\ Soc., to appear.

\bibitem{BrJo97b}
\bysame, \emph{Isometries, shifts, {C}untz algebras and multiresolution wavelet
  analysis of scale ${N}$}, Integral Equations Operator Theory \textbf{28}
  (1997), 382--443.

\bibitem{BJKW97}
O.~Bratteli, P.E.T. Jorgensen, A.~Kishimoto, and R.~Werner, \emph{Pure states
  on $\mathcal{O}_d$}, submitted to Journal of Operator Theory.

\bibitem{BrRoI}
O.~Bratteli and D.W. Robinson, \emph{{O}perator {A}lgebras and {Q}uantum
  {S}tatistical {M}echanics}, 2nd ed., vol.~I, Springer-Verlag, Berlin--New
  York, 1987.

\bibitem{BrRoII}
\bysame, \emph{{O}perator {A}lgebras and {Q}uantum {S}tatistical {M}echanics},
  2nd ed., vol.~II, Springer-Verlag, Berlin--New York, 1996.

\bibitem{CoRy95}
A.~Cohen and R.D. Ryan, \emph{{W}avelets and {M}ultiscale {S}ignal
  {P}rocessing}, Applied Mathematics and Mathematical Computation, vol.~11,
  Chapman \& Hall, London, 1995.

\bibitem{Cun77}
J.~Cuntz, \emph{Simple ${C}^*$-al\-ge\-bras generated by isometries}, Comm.
  Math. Phys. \textbf{57} (1977), 173--185.

\bibitem{DDL95}
Stephan Dahlke, Wolfgang Dahmen, and Vera Latour, \emph{Smooth refinable
  functions and wavelets obtained by convolution products}, Appl. Comput.
  Harmon. Anal. \textbf{2} (1995), 68--84.

\bibitem{Dau92}
I.~Daubechies, \emph{{T}en {L}ectures on {W}avelets}, CBMS-NSF Regional Conf.
  Ser. in Appl. Math., vol.~61, Society for Industrial and Applied Mathematics,
  Philadelphia, 1992.

\bibitem{DaPi97b}
Ken~R. Davidson and David~R. Pitts, \emph{The algebraic structure of
  non-commutative analytic {T}oeplitz algebras}, Math. Ann., to appear.

\bibitem{Dix64}
J.~Dixmier, \emph{Traces sur les ${C}\sp{\ast} $-alg\`ebres, \textup{{II}}},
  Bull. Sci. Math. \textup{(2)} \textbf{88} (1964), 39--57.

\bibitem{Dix69}
\bysame, \emph{Les ${C}^{*}$-alg\`ebres et {L}eurs {R}epr\'esentations},
  Gauthier-Villars, Paris, 1969.

\bibitem{Eva80}
D.E. Evans, \emph{On $\mathcal{O}_{N}$}, Publ. Res. Inst. Math. Sci.
  \textbf{16} (1980), 915--927.

\bibitem{EvKa98}
D.E. Evans and Y.~Kawahigashi, \emph{Quantum {S}ymmetries on {O}perator
  {A}lgebras}, Oxford Mathematical Monographs, The Clarendon Press, Oxford
  University Press, New York, 1998.

\bibitem{GrMa92}
K.~Gr\"ochenig and W.R. Madych, \emph{Multiresolution analysis, {H}aar bases,
  and self-similar tilings of $\mathbf{R}^n$}, IEEE Trans. Inform. Theory
  \textbf{38} (1992), 556--568.

\bibitem{HRW92}
P.N. Heller, H.L. Resnikoff, and R.O. Wells, Jr., \emph{Wavelet matrices and
  the representation of discrete functions}, Wavelets: A Tutorial in Theory and
  Applications (C.K. Chui, ed.), Wavelet Anal. Appl., vol.~2, Academic Press,
  Boston, 1992, pp.~15--50.

\bibitem{HeWe94}
P.N. Heller and R.O. Wells, Jr., \emph{The spectral theory of multiresolution
  operators and applications}, Wavelets: Theory, Algorithms, and Applications
  (Taormina, 1993) (C.K. Chui, L.~Montefusco, and L.~Puccio, eds.), Wavelet
  Anal. Appl., vol.~5, Academic Press, San Diego, 1994, pp.~13--31.

\bibitem{Hor95}
L.~H\"ormander, \emph{Lectures on harmonic analysis}, Dept.\ of Mathematics,
  Box 118, S-22100 Lund, 1995.

\bibitem{JSW94a}
P.E.T. Jorgensen, L.M. Schmitt, and R.F. Werner, \emph{$q$-canonical
  commutation relations and stability of the {C}untz algebra}, Pacific J. Math.
  \textbf{165} (1994), 131--151.

\bibitem{Mal89}
S.G. Mallat, \emph{Multiresolution approximations and wavelet orthonormal bases
  of $\mathbf{L}^{2}(\mathbf{R})$}, Trans. Amer. Math. Soc. \textbf{315}
  (1989), 69--87.

\bibitem{McWe94}
K.~McCormick and R.O. Wells, Jr., \emph{Wavelet calculus and finite difference
  operators}, Math. Comp. \textbf{63} (1994), 155--173.

\bibitem{Pol89}
D.~Pollen, \emph{Parametrization of compactly supported wavelets}, company
  report, AWARE, Inc., AD890503.1.4, 1989.

\bibitem{Pol90}
\bysame, \emph{$\mathrm{SU}\sb{I}(2,{F}[z,1/z])$ for ${F}$ a subfield of
  $\mathbf{C}$}, J.~Amer. Math. Soc. \textbf{3} (1990), 611--624.

\bibitem{Pol92}
\bysame, \emph{{D}aubechies' scaling function on $[0,3]$}, Wavelets: A Tutorial
  in Theory and Applications (C.K. Chui, ed.), Wavelet Anal. Appl., vol.~2,
  Academic Press, Boston, 1992, pp.~3--13.

\bibitem{ReWe92}
H.L. Resnikoff and R.O. Wells, Jr., \emph{Wavelet analysis and the geometry of
  {E}uclidean domains}, J.~Geom. Phys. \textbf{8} (1992), 273--282.

\bibitem{Str94}
R.S. Strichartz, \emph{Self-similarity in harmonic analysis}, J.~Fourier Anal.
  Appl. \textbf{1} (1994), 1--37.

\bibitem{Wel93}
R.O. Wells, Jr., \emph{Parametrizing smooth compactly supported wavelets},
  Trans. Amer. Math. Soc. \textbf{338} (1993), 919--931.

\bibitem{Wic93}
M.V. Wickerhauser, \emph{Best-adapted wavelet packet bases}, Different
  perspectives on wavelets {(San Antonio, TX, 1993)} (I.~Daubechies, ed.),
  Proc. Sympos. Appl. Math., vol.~47, Amer. Math. Soc., Providence, RI, 1993,
  pp.~155--171.

\end{thebibliography}

\end{document}